\documentclass[11pt]{article}
\usepackage{graphicx,amssymb,mathrsfs,amsmath,color,enumerate}
\newfont{\bb}{msbm10}

\newcommand{\tr}{^{\sf T}}

\newtheorem{algorithm}{Algorithm}[section]
\newtheorem{remark}{Remark}[section]

\newtheorem{theorem}{Theorem}[section]
\newtheorem{definition}{Definition}[section]
\newtheorem{lemma}{Lemma}[section]

\begin{document}
\cleardoublepage \pagestyle{plain}
\bibliographystyle{plain}

\title{Accelerated Symmetric ADMM and Its Applications in Signal Processing
\thanks{This research was partially supported by the National Statistical Science Research Project of China (Grant No. 2018LZ23), the Natural Science Foundation of China (Grant No. 11801455, 11571178), Fundamental Research Funds of China West Normal University (Grant No. 17E084, 18B031).}
   }

\author{
Jianchao Bai
\footnote{Department of Applied Mathematics,     Northwestern Polytechnical
University, Xi'an  710129,    China
(\tt jianchaobai@nwpu.edu.cn).}
\quad   Junli Liang$^*$
\footnote{  (Corresponding author) School of Electronics and Information, Northwestern Polytechnical
University, Xi'an  710129,    China
(\tt liangjunli@nwpu.edu.cn).}
\quad Ke Guo
\footnote{    School of Mathematics and Information, China West Normal University, Nanchong, Sichuan 637002,  China
(\tt keguo2014@126.com).}
\quad Yang Jing
\footnote{  School of Electronics and Information, Northwestern Polytechnical
University, Xi'an  710129,    China
(\tt jingyang@mail.nwpu.edu.cn).}
}
\date{June 14th, 2019}
\maketitle

\centerline{\small\it\bf Abstract}\vskip 1mm
The alternating direction method of multipliers (ADMM)    were extensively investigated in the past decades for solving   separable convex optimization problems.  Fewer   researchers focused on exploring  its  convergence  properties for the nonconvex case   although it performed  surprisingly  efficient.  In this paper, we propose a symmetric ADMM based on different acceleration techniques for   a family of potentially nonsmooth nonconvex programing problems with equality constraints, where the dual variables are  updated twice with different stepsizes.  Under proper assumptions instead of  using the so-called Kurdyka-Lojasiewicz inequality, convergence   of the proposed algorithm as well as its   pointwise iteration-complexity are  analyzed     in terms of the corresponding augmented Lagrangian  function and the primal-dual residuals, respectively. Performance of our  algorithm is verified by  some preliminary numerical examples on applications in     sparse nonconvex/convex regularized minimization    signal processing problems.

\vskip 3mm\noindent {\small\bf Keywords:}
Nonconvex optimization,    symmetric ADMM,  acceleration technique,   complexity,  signal processing

\noindent {\small\bf Mathematics Subject Classification(2010):} 47A30; 65Y20;   90C26; 90C90
% 47A30---Norms (inequalities, more than one norm, etc.)
% 65Y20---Complexity and performance of numerical algorithms
% 90C26---Nonconvex programming, global optimization
% 90C90---Applications of mathematical programming
\bigskip

%====================================================================%
\section{Introduction}
%====================================================================%
We consider  a   potentially nonsmooth and nonconvex     separable optimization problem subject to   linear  equality constraints:
\begin{equation} \label{Sec1-Prob}
\min \left\{  f(\mathbf{x})+g(\mathbf{y})|\ \textrm{s.t.} \
    A\mathbf{x}+B\mathbf{y}=b, \mathbf{x}\in \mathcal{R}^m, \mathbf{y}\in \mathcal{R}^n\right\},
\end{equation}
where $f: \mathcal{R}^m\rightarrow (-\infty,+\infty]$ is a  proper lower semicontinuous function, $g: \mathcal{R}^n\rightarrow (-\infty,+\infty)$ is a  continuous differentiable  function  with its gradient $\nabla g$ being $L_g$-Lipschitz continuous, $A\in \mathcal{R}^{l \times m}, B\in \mathcal{R}^{l \times n}, b\in \mathcal{R}^{l}$ are respectively  given matrices and vector.
Minimization problem in the form of (\ref{Sec1-Prob}) covers many important   applications in  science and engineering. For example, the following $l_1$-regularized least square problem arising in signal processing/statistical learning \cite{BoydParikh2010,bAIZ18,KKBG07}:
\begin{equation} \label{example-1}
\min_{x\in\mathcal{R}^m}\frac{1}{2}\|A\mathbf{x}-c\|^2+ \mu\|\mathbf{x} \|_1,
\end{equation}
where  $c\in\mathcal{R}^l$ is the vector of observations, $A\in \mathcal{R}^{l\times m}$ is the data matrix and $\mu>0$ denotes the regularization parameter and is often set     as  $\mu=0.1\mu_{\max}$ where   $\mu_{\max}=\|A\tr c\|_\infty$  (see e.g. \cite{FNW07,KKBG07}). Due to  the convexity of
the   problem (\ref{example-1}), it can be handled by a number of standard   methods,  to list a few, including the   alternating direction method of multipliers (ADMM, \cite{EB92,Gabay83,Glowinski80}),   proximal point algorithm \cite{bAIZ18,EB92},   interior point method \cite{KKBG07} and   primal-dual hybrid gradient method \cite{blw19,zhch08}.  However, in   many cases the $l_1$-regularization has been shown to be sub-optimal. For instance,  it can not recover  a signal with the fewest measurements when being applied in compressed sensing techniques \cite{ChaS08}. Therefore, an  acceptable  improvement is to adopt the $l_{1/2}$-regularization term, which results in the following   form
\[
\min_{x\in\mathcal{R}^m}\frac{1}{2}\|A\mathbf{x}-c\|^2+ \mu\|\mathbf{x} \|_{1/2}^{1/2}.
\]
Here,      $\|\mathbf{x}\|_{1/2}= (\sum_{i=1}^{n}|\mathbf{x}_i|^{\frac{1}{2}})^2$ is a nonconvex function characterizing    sparsity of the variable, and it has been    verified    \cite{xcxz12}  practically to  be better than $l_1$-norm.  Clearly, by introducing an auxiliary variable, the problem can be converted  to a special case of  (\ref{Sec1-Prob}), i.e.,
\begin{equation} \label{example-2}
\min \left\{ \mu\|\mathbf{x} \|_{1/2}^{1/2}+\frac{1}{2}\|\mathbf{y}-c\|^2
|\ \textrm{s.t.} \
    A\mathbf{x}-\mathbf{y}=0  \right\}.
\end{equation}
Another interesting example is the   regularized empirical risk minimization
arising from big data  applications, such as many kinds of classification and regression models
in machine learning \cite{szc13,WLWHam17}. And the  $l_{1/2}$-regularized  reformulation case is of the   form:
\begin{equation} \label{example-3}
\min \left\{\mu\|\mathbf{x}\|_{1/2}^{1/2}+\frac{1}{N}\sum_{j=1}^{N}g_j(\mathbf{y})\\
|\ \textrm{s.t. }   \mathbf{x}-\mathbf{y}=0\right\},
\end{equation}
where   $N$ is a   large number, and $g_j(\mathbf{y})=\log\left(1+\exp(-b_ia_i\tr \mathbf{y})\right)$ denotes the logistic loss function on the feature-label pair $(a_j,b_j)$ with $a_j\in \mathcal{R}^{l}$ and $b_j\in\{-1, 1\}$.

 In the literature, the most standard method for  solving  the equality constrained problem (\ref{Sec1-Prob}) is the augmented Lagrangian method (ALM) which firstly solves  a joint minimization problem
\begin{equation} \label{aug-func}
\min\limits_{\mathbf{x}, \mathbf{y}}\mathcal{L}_\beta(\mathbf{x}, \mathbf{y},\lambda):=f(\mathbf{x})+g(\mathbf{y})-\langle\lambda, A\mathbf{x}+B\mathbf{y}-b\rangle +\frac{\beta}{2}\|A\mathbf{x}+B\mathbf{y}-b\|^2,
\end{equation}
and then updates   the Lagrange multiplier  $\lambda$ by using  the newest iteration of other variables.   The  penalty factor $\beta>0$, in each iterative loop,   can be set as a tuned reasonable value or updated adaptively according to the ratio of the primal residual to the dual residual  of the problem.  However, ALM does not make full use of the separable structure of the objective function of (\ref{Sec1-Prob}) and hence, could not take advantage of the special properties of each component objective function.
This would make it very expensive even infeasible for application problems involving big-data and   nonconvex objectives.  By contrast, a  powerful first-order method, that is ADMM,  aims    to split the joint core  problem (\ref{aug-func}) into some relatively simple and smaller-dimensional subproblems so that variables can be updated separately  to   make full use of  special properties of each component. Another obvious feature of ADMM is that the resultant  subproblems could admit explicit solution form in   special  applications, or in a linearized update for the differentiable objective/quadratic penalty term. We refer to, e.g., \cite{BoydParikh2010,bAIZ18,FGLK15,R14,He-ma-yuan16,HeYangWang2000,WYZ19} for some reviews on ADMM.

Interestingly, under  the existence assumption of a solution to the Karush-Kuhn Tucker condition of the two-block separable convex optimization problem, it was explained
  \cite{Gabay83} that the original ADMM amounts to the Douglas-Rachford
splitting method (DRSM, \cite{DouglasR56,LionsMe79}) when it was  applied to a stationary system to the dual of the problem. Moreover,
 as elaborated  in \cite{Gabay83}, if applying the classic Peaceman-Rachford splitting method (PRSM, \cite{LionsMe79,PeacemanRach55})  to the dual of the problem, we obtain the following  iterative scheme
\begin{equation}  \label{PRSM}
\left \{\begin{array}{lll}
\mathbf{x}_{k+1}&=&\arg\min\limits_{\mathbf{x}} \mathcal{L}_\beta\left(\mathbf{x},\mathbf{y}_k,\lambda_k\right),\\
\lambda_{k+\frac{1}{2}}&=&\lambda_k-  \beta\left(A\mathbf{x}_{k+1}+B\mathbf{y}_{k}-b\right),\\
\mathbf{y}_{k+1}&=&\arg\min\limits_{\mathbf{y}} \mathcal{L}_\beta(\mathbf{x}_{k+1},\mathbf{y},\lambda_{k+\frac{1}{2}}),\\
\lambda_{k+1}&=&\lambda_{k+\frac{1}{2}}- \beta\left(A\mathbf{x}_{k+1}+B\mathbf{y}_{k+1}-b\right).
\end{array}\right.
\end{equation}
Unfortunately, scheme \eqref{PRSM} is not convergent under the standard  convexity assumptions as  ADMM \cite{CY14}. However,  it was verified  \cite{GlowinskiKMa03} that scheme \eqref{PRSM}  could perform faster than the ADMM when  its global  convergent was ensured.   In view of this,  He et al. in \cite{HLWY14} proposed and studied the convergence of a strictly contractive Peaceman-Rachford splitting method (also called the symmetric version of ADMM)
\begin{equation}  \label{SPRSM}
\left \{\begin{array}{lll}
\mathbf{x}_{k+1}&=&\arg\min\limits_{\mathbf{x}} \mathcal{L}_\beta\left(\mathbf{x},\mathbf{y}_k,\lambda_k\right),\\
\lambda_{k+\frac{1}{2}}&=&\lambda_k-  \alpha\beta\left(A\mathbf{x}_{k+1}+B\mathbf{y}_{k}-b\right),\\
\mathbf{y}_{k+1}&=&\arg\min\limits_{\mathbf{y}} \mathcal{L}_\beta(\mathbf{x}_{k+1},\mathbf{y},\lambda_{k+\frac{1}{2}}),\\
\lambda_{k+1}&=&\lambda_{k+\frac{1}{2}}- \alpha\beta\left(A\mathbf{x}_{k+1}+B\mathbf{y}_{k+1}-b\right),
\end{array}\right.
\end{equation}
where $\alpha\in (0,1)$ is the relaxation parameter.  Later, He et al. \cite{He-ma-yuan16} improved the scheme (\ref{SPRSM})   to the case  with  larger range of relaxation parameters, which was   generalized  by Bai et al. \cite{BaiZhangLiXu2016}  to the multi-block separable convex programming. Besides,   Chang, et. al.\cite{ChangL19} also shown a  generalization of linearized ADMM for  two-block separable convex minimization model by adding a proper  proximal term to each core subproblem.

If the convexity is lose, then the convergence analysis for ADMM (or its variant) is much more challenging.  However, for some special nonconvex optimization problems, one can establish  convergence of   ADMM by making full use of special structures of the  problems,   see e.g. \cite{HLR16} for the consensus and sharing problems.   Another widely used technique to prove   convergence of   ADMM for nonconvex optimization problems relies on the assumption that the objective function of (1) satisfies the so-called Kurdyka-Lojasiewicz (KL) inequality \cite{Attbrs10},   since many important classes of functions satisfy the KL inequality,   see \cite{GHW17,GHWW17,GHW18,LiPONG15,WYZ19,WLWHam17,YPC17}.
 Without assuming the KL property and convexity  of the objective function, recently,  Goncalves et al. in  \cite{GMM17} established convergence rate bounds of the classical ADMM with  proximal terms for solving nonconvex linearly constrained optimization problem  (\ref{Sec1-Prob}). In addition, by linearizing the smooth part in the objective and  quadratic penalty term, Liu, et al. \cite{LSGU17} proposed a two-block linearized ADMM for the problem (\ref{Sec1-Prob}) with $b=0$ and extended  the method   to a multi-block version, but convergence of their  extended  method holds with an extra  hypothesis on the full column rank of  the matrix $B$ compared to (A1) (see Section 3).

Motivated by the above mentioned work  \cite{GMM17,LSGU17} and the empirical validity of the symmetric ADMM, we would present a  Two-stage Accelerated Symmetric ADMM (abbreviated as ``TAS-ADM") for solving the problem (\ref{Sec1-Prob}), whose framework reads    Algorithm \ref{algo1}.  Our algorithm combines both the so-called Nesterov's acceleration technique  in (\ref{jb--2})  and the relaxation  scheme  in e.g., \cite{EB92,FGLK15}.  By adding a proper proximal term for the first $\mathbf{x}$-subproblem, this possibly nonsmooth nonconvex subproblem will turn to a proximal mapping   shown in (\ref{bj-40}), which admits closed solution form if $f$ is easy.  Step 7 actually uses the idea of  convex combination for fast convergence.

We should emphasize that the recent work  \cite{WLWHam17} also considered a symmetric ADMM  for solving the problem (\ref{Sec1-Prob}). The method in \cite{WLWHam17} actually can be treated as our proposed Algorithm \ref{algo1} barring the acceleration techniques and proximal regularization terms, while convergence of Algorithm \ref{algo1} is analyzed in a different way.  More precisely, their analyses  are based on the Kurdyka-Lojasiewicz property of the augmented Lagrangian function for problem (\ref{Sec1-Prob}) and other proper assumptions on both the penalty parameter and the objective function.
 Under Assumptions (A1-A3) (see Section 3), we show in the sequel section  that any accumulation point of $\{w_k:=(\mathbf{x}_k,\mathbf{y}_k,\lambda_k)\}$ is the stationary point of $\{L_\beta(w_k)\}$, and we also establish the worst-case $\mathcal{O}(1/k)$ convergence rate of the algorithm in terms of the primal-dual residuals.   Although we consider problem (\ref{Sec1-Prob}) with vector variables,  the subsequent  convergence results of our proposed algorithm are applicable for the general case with matrix variables, because matrix can be vectorized as vector.

%By adding a proper positive (in)definite proximal term to each core subproblem and updating the Lagrangian multiplier twice in different ways,  Chang, et. al.\cite{ChangL19} also shown a  generalization of linearized ADMM for solving two-block separable convex minimization model.

\begin{algorithm}\label{algo1}
\vskip1mm
\hrule\vskip2mm
[TAS-ADM for Solving Problem (\ref{Sec1-Prob})]\vskip1.3mm
\hrule\vskip1.3mm
\noindent \verb"1  Initialize"   $(\mathbf{x}_0,\mathbf{y}_0,\lambda_0)\in \mathcal{R}^m\times \mathcal{R}^n\times \mathcal{R}^l$ \verb"and set" $(\mathbf{x}_{-1},\mathbf{y}_{-1})=(\mathbf{x}_0,\mathbf{y}_0).$ \\
\verb"2  Choose parameters" $\beta>0, \gamma_k\in[0,\frac{1}{2})$,   $G \succeq 0 $ \verb"and" \\
\begin{equation} \label{al-111}
(\tau, \alpha )\in \mathcal{D}:=\left\{(\tau, \alpha )|\ 0<\tau+\alpha <1\right\}.
\end{equation}
\verb"3  for"  $k=0,1,\cdots,$   \verb"do"  \\
\verb"4"\indent\quad\ \
$\mathbf{x}_{k}^{md}=\mathbf{x}_{k}+\gamma_{k}(\mathbf{x}_{k}-\mathbf{x}_{k-1}).$\\
\verb"5" \indent\quad\
 $\mathbf{x}_{k+1}  =\arg\min \left\{\mathcal{L}_\beta(\mathbf{x}, \mathbf{y}_k,\lambda_k)+\frac{1}{2}\|\mathbf{x}-\mathbf{x}_{k}^{md}\|_{G}^2\right\}.$ \\
\verb"6" \indent\quad\
$\lambda_{k+\frac{1}{2}}=\lambda_k- \tau \beta\left(A\mathbf{x}_{k+1}+B\mathbf{y}_{k}-b\right).$\\
\verb"7" \indent\quad\
$\mathbf{x}_{k+1}^{ad}  =\alpha  A\mathbf{x}_{k+1}+(1-\alpha )(b-B\mathbf{y}_k).$\\
\verb"8" \indent\quad\
$\mathbf{y}_{k+1}=\arg\min \left\{g(\mathbf{y})-\left\langle\lambda_{k+\frac{1}{2}}, B\mathbf{y}\right\rangle+\frac{\beta}{2}\left\|\mathbf{x}_{k+1}^{ad}
+B\mathbf{y}-b\right\|^2\right\}.$\\
\verb"9" \indent\quad\
$\lambda_{k+1}=\lambda_{k+\frac{1}{2}}-  \beta\left(\mathbf{x}_{k+1}^{ad}+B\mathbf{y}_{k+1}-b\right).$\\
\verb"10 end"\\
\verb"11 Output" $(\mathbf{x}_{k+1},\mathbf{y}_{k+1})$.
\vskip2mm\hrule\vskip3mm
\end{algorithm}

The remaining parts of this paper are organized as follows. In Section \ref{Sepre},  some preliminaries are prepared  to analyze convergence of  Algorithm \ref{algo1}. In Section \ref{Theor res}, we show its convergence properties  and its pointwise iteration complexity based on   the analysis for the augmented Lagrangian sequence   $\{L_\beta(w_k)\}$. Section \ref{numex} tests  some examples about the popular sparse signal recovery problem with different regularization terms and compared with the popular CVX toolbox, which aims to investigate  numerical performance of our   algorithm. Finally, we   conclude the paper in Section \ref{Fin-5}.

%==============================================================================
\section{Preliminaries}\label{Sepre}
%==============================================================================
Throughout this paper, let  $\mathcal{R}, \mathcal{R}^n,\mathcal{R}^{m\times n}$
be the sets of real numbers,  $n$ dimensional real column vectors and
 $m\times n$ dimensional real matrices, respectively.  The symbol  $I$ denotes the identity matrix  with proper
dimension and   $\sigma_B$ denotes the smallest positive eigenvalue of the matrix $BB\tr$.    For any symmetric matrices $A$ and $B$ whose dimensions are the same, $A \succ B$
($A \succeq B$) means $A - B$ is a positive definite (semidefinite) matrix.
We slightly denote $\|x\|_G^2 = x \tr G x$ for any symmetric matrix $G$, and
let $\|x\|_G = \sqrt{x \tr G x}$ when $G$ is positive semidefinite,
where the superscript   $\tr$ denotes the transpose  of a matrix or vector. We simply  use $\|\cdot\|$ to represent
 the standard Euclidean norm  equipped with inner product $\langle \cdot,\cdot\rangle$.
    The image space of a matrix $A\in\mathcal{R}^{m\times n}$ is defined as
 $\textrm{Im}(A):=\{Az|\ z\in \mathcal{R}^n\}$ and a function $f:\mathcal{S}\rightarrow \mathcal{R}$ is lower semicontinuous at $\bar{x}\in \mathcal{S}$ if and only if $\liminf\limits_{x\rightarrow \bar{x}} f(x)=f(\bar{x}).$  The distance from any point $z$ to the set $\mathcal{S}\subseteq \mathcal{R}^n$ is defined as $d(z,\mathcal{S}):=\inf\{\|z-y\|\ | y\in \mathcal{S}\}.$

\begin{definition}\label{def-0}\cite{Mord6, Roc72}
Let $f: \mathcal{R}^m\rightarrow \mathcal{R}$ be a  proper lower semicontinuous function.
\begin{enumerate}[{~~~~(}a{)}]
\item
For a given $x\in \textrm{dom} (f)$, the Frechet subdifferential of $f$ at $x$, written by   $\widehat{\partial} f(x)$, is the set of all vectors $s\in \mathcal{R}^m$  which satisfy
\[
\lim_{y\neq x,}\inf_{y\rightarrow x} \frac{f(y)-f(x)-\langle s, y-x\rangle}{\|y-x\|}\geq 0,
\]
and we let  $\widehat{\partial} f(x)=\emptyset$ when $x \notin\textrm{dom} (f)$.
\item
The limiting subdifferential, or   the subdifferential of $f$ at $x\in \mathcal{R}^m$, written by $\partial f(x)$, is defined by
$
\partial f(x)=\left\{s\in\mathcal{R}^m|\ \exists x_k\rightarrow x, f(x_k)\rightarrow f(x), \widehat{\partial} f(x_k)\ni s_k \rightarrow s  \ \textrm{as}  \ k\rightarrow\infty\right\}.
$
\item
A  point $x_*$ is called  critical point or  stationary point of $f(x)$ if it satisfies $0\in \partial f(x_*)$.
\end{enumerate}
\end{definition}

\begin{definition}\label{def-1}
A triple $w_*:=(\mathbf{x}_*,\mathbf{y}_*,\lambda_*)\in \mathcal{R}^m\times \mathcal{R}^n \times \mathcal{R}^l$ is a stationary point of   (\ref{Sec1-Prob}) if
\[
A\tr \lambda_*\in \partial f(\mathbf{x}_*),\quad
B\tr \lambda_*=\nabla g(\mathbf{y}_*) \quad  \textrm{and} \quad
A\mathbf{x}_*+B\mathbf{y}_*-b=0.
\]
\end{definition}

The following     lemmas are    provided
to simplify  convergence   analysis in the sequel sections.
\begin{lemma}\label{lem-3} \cite[Lemma A.2]{GMM17}
Let $A\in \mathcal{R}^{m\times n}$ be a  nonzero matrix and $\mathcal{P}_{A}$ be  the Euclidean projection onto $\textrm{Im}(A).$
Then, for any $u\in \mathcal{R}^n$ we have
\begin{equation}\label{lem-3-02}
\|\mathcal{P}_{A}(u)\|\leq \frac{1}{\sqrt{\sigma_A}}\|A\tr u\|.
\end{equation}
\end{lemma}

\begin{lemma}\label{lem-1}
For any  vectors  $a,b,c\in \mathcal{R}^n$ and   symmetric matrix $0\preceq M\in \mathcal{R}^{n\times n}$, it  holds
\begin{equation}\label{lem-1-1}
\langle a-b, M(a-c)\rangle =\frac{1}{2}\left\{ \|c-a\|_M^2 - \|c-b\|_M^2 +\|a-b\|_M^2\right\}.
\end{equation}
\end{lemma}
%============================================================================
\section{Theoretical Results}\label{Theor res}
%============================================================================
In this section, by making use  of the following primal-dual residuals
 \begin{equation} \label{aug-func-new}
\triangle\mathbf{x}_k= \mathbf{x}_{k}-\mathbf{x}_{k-1},\quad \triangle\mathbf{y}_k= \mathbf{y}_{k}-\mathbf{y}_{k-1}\quad \textrm{and}\quad \triangle\lambda_k= \lambda_{k}-\lambda_{k-1},
 \end{equation}
 the  proposed algorithm will be demonstrated  to be convergent according to a quasi-monotonically nonincreasing property of the sequence $\{\mathcal{L}_\beta(w_k)\}$, and its  pointwise iteration-complexity will be established in detail.    Next, we make some assumptions.
\begin{itemize}
\item (\textbf{A1})~
$B\neq 0, \ \textrm{Im}(B)\supset {b}\cup \textrm{Im}(A)$;
\item
(\textbf{A2})~The penalty parameter $\beta$ satisfies
\begin{align*}
\beta>\frac{L_{g}}{\sqrt{1-\tau-\alpha}\sigma_{B}},~~~(\tau,\alpha)\in \mathcal{D}~ \textrm{with}~ \mathcal{D}~ \textrm{given\ in~} (\ref{al-111});
\end{align*}

\item (\textbf{A3})~
$ \underline{L}=\inf\limits_{(\mathbf{x},\mathbf{y})}\left\{ f(\mathbf{x})+ g(\mathbf{y})-\frac{1}{2L_g}\|\nabla g(\mathbf{y})\|^2\right\}>-\infty.$
\end{itemize}
Indeed, we can check that the aforementioned Assumptions (A1)-(A3) hold for the two examples mentioned in the introduction.
 Here and hereafter, we   denote
$
w_k=(\mathbf{x}_k,\mathbf{y}_k,\lambda_k)$ and $     w=(\mathbf{x},\mathbf{y},\lambda).
$

\begin{lemma} \label{bj-1}
Let  $\{w_k\}$ be generated by Algorithm \ref{algo1}. Then, under (A2) we have
 \begin{equation} \label{bj-2}
\|B\tr \triangle\lambda_{k+1}\|\leq L_g \|\triangle\mathbf{y}_{k+1}\|.
\end{equation}
\end{lemma}
\noindent{\bf Proof }
According to the optimality condition of $\mathbf{y}$-subproblem, it holds
 \begin{equation} \label{bj-002}
\nabla g(\mathbf{y}_{k+1}) -B\tr\lambda_{k+\frac{1}{2}} +\beta B\tr \left(\mathbf{x}_{k+1}^{ad}+B\mathbf{y}_{k+1}-b\right)=0.
\end{equation}
So, we have by the update  of $\lambda_{k+1}$   that
\begin{equation}\label{bj-3}
B\tr\lambda_{k+1} =  B\tr\left[\lambda_{k+\frac{1}{2}}-\beta\left(\mathbf{x}_{k+1}^{ad}+B\mathbf{y}_{k+1}-b\right)\right]=\nabla g(\mathbf{y}_{k+1}),
\end{equation}
which further gives
\begin{equation} \label{bj-4}
B\tr\lambda_{k} =\nabla g(\mathbf{y}_{k}).
\end{equation}
Subtracting (\ref{bj-4}) from (\ref{bj-3})  and taking norm on both sides, we can obtain  by (A2) that
\begin{eqnarray*}
\|B\tr \triangle\lambda_{k+1}\|
\leq  L_g \|\triangle\mathbf{y}_{k+1}\|. \ \ \ \blacksquare
\end{eqnarray*}

Note that   optimality condition of the following problem is the same as (\ref{bj-002}):
\[
\min \left\{g(\mathbf{y})+\frac{\beta}{2}\left\|B\mathbf{y}-c_y\right\|^2\right\},
\]
where $c_y=b+\frac{\lambda_{k+\frac{1}{2}}}{\beta} - \mathbf{x}_{k+1}^{ad}$. So, this problem is equivalent to the $\mathbf{y}$-subproblem in Algorithm \ref{algo1}. Under the case that $g$ is linearized or $B$ has full column rank, the above problem could  have  closed solution form.
In addition, by choosing $G=\sigma I -\beta A\tr A$ with $\sigma\geq \beta\|A\tr A\|$,  the quadratic term $\|A\mathbf{x}\|^2$ will be cancelled in the iteration. As a result,
 the $\mathbf{x}$-subproblem in Algorithm \ref{algo1} is converted to a  proximal mapping as the following
\begin{equation} \label{bj-40}
\textrm{Prox}_{f,\sigma}(c_x):=\textrm{Arg}\min \left\{f(\mathbf{x})+\frac{\sigma}{2}\left\|\mathbf{x}-c_x\right\|^2\right\},
\end{equation}
where $c_x=\mathbf{x}_k^{md}-\frac{  \beta A\tr(A\mathbf{x}_k^{md} +B\mathbf{y}_k-b)-A\tr\lambda_k  }{\sigma}.$ Since $f$ is a  proper lower semicontinuous function  and bounded from below (in view of Assumptions (A3)), by the proximal behavior in \cite{Roc72} the set $\textrm{Prox}_{f,\sigma}(c_x)$  is  nonempty and compact.

Now, adding the update of $\lambda_{k+\frac{1}{2}}$ to the update of $\lambda_{k+1}$, we have
\begin{eqnarray*}
\frac{1}{\beta}\triangle \lambda_{k+1}
&=& -\tau (A\mathbf{x}_{k+1}+B \mathbf{y}_{k}-b) - \left[\alpha  A\mathbf{x}_{k+1}+(1-\alpha )(b-B\mathbf{y}_k)+B\mathbf{y}_{k+1}-b\right]\\
&=& -(\tau+\alpha )(A\mathbf{x}_{k+1}+B \mathbf{y}_{k}-b)-B\triangle\mathbf{y}_{k+1},
\end{eqnarray*}
which by  $\tau+\alpha >0$ gives the following lemma   immediately.

\begin{lemma} \label{bj-1-01}
Assume $\tau+\alpha >0$, then the sequence $\{w_k\}$   generated by Algorithm \ref{algo1} satisfies
 \begin{equation} \label{bj-1-02}
A\mathbf{x}_{k+1}+B\mathbf{y}_k-b=-\frac{1}{\tau+\alpha }\left(\frac{1}{\beta}\triangle\lambda_{k+1}+B\triangle\mathbf{y}_{k+1}\right).
\end{equation}
\end{lemma}

Next, we present   a fundamental lemma that plays a key role in analyzing convergence and convergence rate bound of Algorithm \ref{algo1}.

\begin{lemma} \label{Sec3-theore1}
Under  Assumptions (A1) and (A2), there exist  three constants
 $\zeta_0\geq0$ and $ \zeta_1,\zeta_2> 0$ such that
\begin{eqnarray} \label{bj-5}
\widetilde{L}_\beta(w_k)- \widetilde{L}_\beta(w_{k+1})  \geq \zeta_1\|\triangle\mathbf{x}_{k+1}\|_G^2+
 \zeta_2 \|\triangle \mathbf{y}_{k+1}\|^2,
\end{eqnarray}
where $\widetilde{L}_\beta(w_k):=\mathcal{L}_\beta(w_k)+ \zeta_0\|\triangle\mathbf{x}_{k}\|_G^2.$
\end{lemma}
\noindent{\bf Proof } The inequality (\ref{bj-5}) can be proved  by the following four steps.

 (\textbf{Step 1}) By the update of $\mathbf{x}$-subproblem together with the way of generating $\mathbf{x}_{k}^{md}$, we have
\begin{eqnarray} \label{bj-6}
&&\mathcal{L}_\beta(\mathbf{x}_{k}, \mathbf{y}_{k},\lambda_{k})- \mathcal{L}_\beta(\mathbf{x}_{k+1}, \mathbf{y}_{k},\lambda_{k})\nonumber\\
&\geq& \frac{1}{2}\left[\|\mathbf{x}_{k+1}-\mathbf{x}_{k}^{md}\|_G^2-\| \mathbf{x}_{k}-\mathbf{x}_{k}^{md}\|_G^2\right]\nonumber\\
&=&\frac{1}{2}\left[\|\mathbf{x}_{k+1}-\mathbf{x}_{k}\|_G^2+2\left\langle \mathbf{x}_{k+1}-\mathbf{x}_{k}, G(\mathbf{x}_{k}-\mathbf{x}_{k}^{md})\right\rangle\right]\nonumber\\
&=&\frac{1}{2}\left[\|\triangle\mathbf{x}_{k+1}\|_G^2-2\gamma_k\left\langle \triangle\mathbf{x}_{k+1}, G\triangle\mathbf{x}_{k}\right\rangle\right]\nonumber\\
&\geq& \frac{1}{2}\left[\|\triangle\mathbf{x}_{k+1}\|_G^2- \gamma_k\left( \|\triangle\mathbf{x}_{k+1}\|_G^2+\|\triangle\mathbf{x}_{k}\|_G^2\right)\right]\nonumber\\
&=&\zeta_0\left[\|\triangle\mathbf{x}_{k+1}\|_G^2-\|\triangle\mathbf{x}_{k}\|_G^2 \right] + \zeta_1\|\triangle\mathbf{x}_{k+1}\|_G^2,
\end{eqnarray}
where
\begin{equation} \label{bj-6-01}
\zeta_0=\frac{\gamma_k}{2}\geq 0,\quad \textrm{and} \quad  \zeta_1= \frac{1-2\gamma_k}{2}> 0.
\end{equation}

(\textbf{Step 2}) By the update of $\mathbf{y}$-subproblem we obtain
\[
g(\mathbf{y}_{k})-\langle \lambda_{k+\frac{1}{2}}, B\mathbf{y}_{k}\rangle + \frac{\beta}{2}\|\mathbf{x}_{k+1}^{ad}+B\mathbf{y}_{k}-b\|^2 \geq
g(\mathbf{y}_{k+1})-\langle \lambda_{k+\frac{1}{2}}, B\mathbf{y}_{k+1}\rangle +\frac{\beta}{2}\|\mathbf{x}_{k+1}^{ad}+B\mathbf{y}_{k+1}-b\|^2,
\]
which, by   Lemma \ref{lem-1}, is equivalently expressed as
\begin{equation} \label{bj-6-1}
 g(\mathbf{y}_{k})-g(\mathbf{y}_{k+1})+ \langle \lambda_{k+\frac{1}{2}}, B\triangle\mathbf{y}_{k+1}\rangle+ \frac{\beta}{2}\|B\triangle\mathbf{y}_{k+1}\|^2
 \geq \beta\langle B\triangle\mathbf{y}_{k+1}, \mathbf{x}_{k+1}^{ad}+B\mathbf{y}_{k+1}-b\rangle.
\end{equation}
Therefore, it  can be deduced that
\begin{eqnarray} \label{bj-7}
&&\mathcal{L}_\beta(\mathbf{x}_{k+1}, \mathbf{y}_{k},\lambda_{k+\frac{1}{2}})- \mathcal{L}_\beta(\mathbf{x}_{k+1}, \mathbf{y}_{k+1},\lambda_{k+\frac{1}{2}})\nonumber\\
&=& g(\mathbf{y}_{k})- g(\mathbf{y}_{k+1})+ \left\langle \lambda_{k+\frac{1}{2}},B\triangle\mathbf{y}_{k+1}\right\rangle+ \frac{\beta}{2}\left(\|A\mathbf{x}_{k+1}+B\mathbf{y}_{k}-b\|^2-\|A\mathbf{x}_{k+1}+B\mathbf{y}_{k+1}-b\|^2\right)\nonumber\\
&=& g(\mathbf{y}_{k})- g(\mathbf{y}_{k+1})+ \left\langle \lambda_{k+\frac{1}{2}},B\triangle\mathbf{y}_{k+1}\right\rangle
-\beta\langle B\triangle\mathbf{y}_{k+1}, A\mathbf{x}_{k+1}+B\mathbf{y}_{k+1}-b\rangle+ \frac{\beta}{2}\|B\triangle\mathbf{y}_{k+1}\|^2\nonumber\\
&\geq & \beta\langle B\triangle\mathbf{y}_{k+1}, \mathbf{x}_{k+1}^{ad}+B\mathbf{y}_{k+1}-b\rangle -\beta\langle B\triangle\mathbf{y}_{k+1}, A\mathbf{x}_{k+1}+B\mathbf{y}_{k+1}-b\rangle\nonumber\\
&=&  \beta(\alpha  -1)\langle B\triangle\mathbf{y}_{k+1}, A\mathbf{x}_{k+1}+ B\mathbf{y}_{k}-b\rangle \nonumber\\
&=& \frac{1-\alpha }{\tau+\alpha }\left[\beta\|B\triangle\mathbf{y}_{k+1}\|^2+ \langle\triangle\mathbf{y}_{k+1}, B\tr\triangle\lambda_{k+1}\rangle\right],
\end{eqnarray}
where the second equality follows   Lemma \ref{lem-1}, the first inequality uses (\ref{bj-6-1}), the third equality uses the update of $\mathbf{x}_{k+1}^{ad}$ and the final equality uses (\ref{bj-1-02}).

(\textbf{Step 3}) Note that
\begin{eqnarray} \label{bj-8}
&&  \mathcal{L}_\beta(\mathbf{x}_{k+1}, \mathbf{y}_{k},\lambda_{k})-  \mathcal{L}_\beta(\mathbf{x}_{k+1}, \mathbf{y}_{k},\lambda_{k+\frac{1}{2}})+
\mathcal{L}_\beta(\mathbf{x}_{k+1}, \mathbf{y}_{k+1},\lambda_{k+\frac{1}{2}})- \mathcal{L}_\beta(\mathbf{x}_{k+1}, \mathbf{y}_{k+1},\lambda_{k+1})\nonumber\\
&=&
 \left\langle \lambda_{k+\frac{1}{2}}-\lambda_{k}, A\mathbf{x}_{k+1}+B\mathbf{y}_{k}-b \right\rangle - \left\langle \lambda_{k+\frac{1}{2}}-\lambda_{k+1}, A\mathbf{x}_{k+1}+B\mathbf{y}_{k+1}-b  \right\rangle\nonumber\\
&=&  \left\langle \lambda_{k+\frac{1}{2}}-\lambda_{k}, A\mathbf{x}_{k+1}+B\mathbf{y}_{k}-b \right\rangle - \left\langle \lambda_{k+\frac{1}{2}}\underbrace{- \lambda_{k}+ \lambda_{k}}-\lambda_{k+1}, A\mathbf{x}_{k+1}+B\mathbf{y}_{k+1}-b  \right\rangle\nonumber\\
&=& \left\langle \lambda_{k+\frac{1}{2}}-\lambda_{k}, -B\triangle \mathbf{y}_{k+1}\right\rangle + \left\langle \triangle\lambda_{k+1}, A\mathbf{x}_{k+1}+B\mathbf{y}_{k+1}-b  \right\rangle\nonumber\\
&=& \tau\beta\left\langle A\mathbf{x}_{k+1}+B \mathbf{y}_{k}-b, B\triangle \mathbf{y}_{k+1}\right\rangle + \left\langle \triangle\lambda_{k+1}, A\mathbf{x}_{k+1}+B\mathbf{y}_{k+1}\underbrace{+B\mathbf{y}_{k}-B\mathbf{y}_{k}}-b  \right\rangle\nonumber\\
&=& \left\langle A\mathbf{x}_{k+1}+B \mathbf{y}_{k}-b, \triangle\lambda_{k+1}+\tau\beta B\triangle \mathbf{y}_{k+1}\right\rangle +
\left\langle \triangle\lambda_{k+1}, B\triangle \mathbf{y}_{k+1}\right\rangle \nonumber\\
&=& -\frac{1}{\tau+\alpha }\left\langle \frac{1}{\beta}\triangle \lambda_{k+1} +B\triangle \mathbf{y}_{k+1}, \triangle\lambda_{k+1}+\tau\beta B\triangle \mathbf{y}_{k+1}\right\rangle +
\left\langle \triangle\lambda_{k+1}, B\triangle \mathbf{y}_{k+1}\right\rangle \nonumber\\
&=&  -\frac{\tau  \beta}{\tau+\alpha } \|B\triangle \mathbf{y}_{k+1}\|^2 -\frac{1}{(\tau+\alpha )\beta} \left\|\triangle \lambda_{k+1}\right\|^2-\frac{1-\alpha }{\tau+\alpha }\left\langle   \triangle\lambda_{k+1}, B\triangle \mathbf{y}_{k+1}\right\rangle
.
\end{eqnarray}

(\textbf{Step 4})
Summing the above inequalities (\ref{bj-6}), (\ref{bj-7}) and the equality (\ref{bj-8}), we  get
\begin{eqnarray*}
&&\mathcal{L}_\beta(\mathbf{x}_{k}, \mathbf{y}_{k},\lambda_{k})-\mathcal{L}_\beta(\mathbf{x}_{k+1}, \mathbf{y}_{k+1},\lambda_{k+1})\\
 & \geq&  \frac{\gamma_k}{2}\left[\|\triangle\mathbf{x}_{k+1}\|_G^2-\|\triangle\mathbf{x}_{k}\|_G^2 \right] + \frac{1-2\gamma_k}{2}\|\triangle\mathbf{x}_{k+1}\|_G^2 +R_\triangle,
\end{eqnarray*}
where
\begin{eqnarray*}
R_\triangle&=&\left(\frac{1}{\tau+\alpha }-1 \right)\beta\|B\triangle \mathbf{y}_{k+1}\|^2
   -\frac{1}{(\tau+\alpha )\beta} \left\|\triangle \lambda_{k+1}\right\|^2\\
   &\geq&  \left(\frac{1}{\tau+\alpha }-1 \right)\beta\|B\triangle \mathbf{y}_{k+1}\|^2
   -\frac{1}{(\tau+\alpha )\beta\sigma_B} \left\|B\tr\triangle \lambda_{k+1}\right\|^2\\
   &\geq&  \left(\frac{1}{\tau+\alpha }-1 \right)\beta\|B\triangle \mathbf{y}_{k+1}\|^2
   -\frac{L_g^2}{(\tau+\alpha )\beta\sigma_B}  \|\triangle \mathbf{y}_{k+1}\|^2\\
   &\geq&  \left(\frac{1}{\tau+\alpha }-1 \right)\beta\sigma_B\|\triangle \mathbf{y}_{k+1}\|^2
   -\frac{L_g^2}{(\tau+\alpha )\beta\sigma_B}  \|\triangle \mathbf{y}_{k+1}\|^2\\
   &=& \zeta_2\|\triangle \mathbf{y}_{k+1}\|^2
\end{eqnarray*}
with
\[
 \zeta_2=\frac{(1-\tau-\alpha )\beta^2\sigma_B^2-L_g^2}{(\tau+\alpha )\beta\sigma_B}> 0. \textrm{~ {[due\ to\ (A2)]}}
\]
Actually, in the first inequality of   $R_\triangle$,
we use the fact that $ \triangle \lambda_{k+1}\in \textrm{Im}(B)$ because of   Assumption (A1).
So,  the  whole proof is completed by the notation $\widetilde{L}_\beta(w_k)$.
 $\ \ \ \blacksquare$

\begin{theorem} \label{lem-34}
Let
 $\{w_k\}$ be generated by Algorithm \ref{algo1}. Then, under   (A1)-(A3) we have
 \begin{itemize}
\item
  The sequence $\{L_\beta(w_k)\}$ is convergent;
\item
The residuals $\|\triangle\mathbf{x}_{k+1}\|_G$, $\|\triangle\mathbf{y}_{k+1}\|$ \textrm{and}
$\|\triangle\lambda_{k+1}\|$ converge  to zero as $k$ goes to infinity.
\end{itemize}
\end{theorem}
\noindent{\bf Proof } To demonstrate   convergence  of $\{L_\beta(w_k)\}$, we need to make ensure  that
  the sequence $\{w_k\}$ is bounded at first.
  By Assumption (A2), it holds
\[
L_g<\sqrt{1-\tau-\alpha }\beta \sigma_B< \beta \sigma_B.
\]
Combining the above inequality and   Lemma \ref{Sec3-theore1},  we achieve
\begin{eqnarray} \label{lema31-01}
&&L_\beta(\mathbf{x}_0,\mathbf{y}_0,\lambda_0)
=  L_\beta(\mathbf{x}_0,\mathbf{y}_0,\lambda_0)+ \zeta_1\|\triangle\mathbf{x}_{0}\|_G^2 \nonumber\\
&\geq&  L_\beta(\mathbf{x}_{k+1},\mathbf{y}_{k+1},\lambda_{k+1})+  \zeta_1\|\triangle\mathbf{x}_{k+1}\|_G^2\geq  L_\beta(\mathbf{x}_{k+1},\mathbf{y}_{k+1},\lambda_{k+1}) \nonumber\\
&=& f(\mathbf{x}_{k+1})+g(\mathbf{y}_{k+1})-\frac{1}{2\beta}\|\lambda_{k+1}\|^2 +\frac{\beta}{2}\left\|A\mathbf{x}_{k+1}+B\mathbf{y}_{k+1}-b
-\frac{\lambda_{k+1}}{\beta}\right\|^2 \nonumber\\
&\geq& f(\mathbf{x}_{k+1})+g(\mathbf{y}_{k+1})-\frac{1}{2\beta\sigma_B}\|B\tr\lambda_{k+1}\|^2 +\frac{\beta}{2}\left\|A\mathbf{x}_{k+1}+B\mathbf{y}_{k+1}-b
-\frac{\lambda_{k+1}}{\beta}\right\|^2 \nonumber\\
&=&\left(f(\mathbf{x}_{k+1})+g(\mathbf{y}_{k+1})-\frac{1}{2L_g}\|\nabla g(\mathbf{y}_{k+1})\|^2\right) +\left(\frac{1}{2L_g}
-\frac{1}{2\beta\sigma_B}\right)\|B\tr\lambda_{k+1}\|^2 \nonumber\\
&&+\frac{\beta}{2}\left\|A\mathbf{x}_{k+1}+B\mathbf{y}_{k+1}-b
-\frac{\lambda_{k+1}}{\beta}\right\|^2 \nonumber\\
&\geq&  \underline{L}   +\left(\frac{1}{2L_g}
-\frac{1}{2\beta\sigma_B}\right)\|B\tr\lambda_{k+1}\|^2 +\frac{\beta}{2}\left\|A\mathbf{x}_{k+1}+B\mathbf{y}_{k+1}-b
-\frac{\lambda_{k+1}}{\beta}\right\|^2,
\end{eqnarray}
which implies that the sequences $\{\lambda_{k}\}, \{\frac{\beta}{2}\left\|A\mathbf{x}_{k+1}+B\mathbf{y}_{k+1}-b
- {\lambda_{k+1}}/{\beta}\right\|^2\}$ are bounded, and furthermore  both $\{\mathbf{x}_{k}\}$ and $\{\mathbf{y}_{k}\}$ are bounded. So, the sequence $\{w_k\}$ is bounded.

Since $\{w_k\}$ is bounded,   $\{\widetilde{L}_\beta(w_k)\}$ is also bounded from below and there exists at least one limit point. Without loss of generality, let $w_*$ be the limit point of $\{w_k\}$ whose subsequence is   $\{w_{k_j}\}$. Then,  the lower semicontinuity of   $\{\widetilde{L}_\beta(w)\}$ indicates
\[
\widetilde{L}_\beta(w_*)\leq \liminf\limits_{j\rightarrow +\infty}\widetilde{L}_\beta(w_{k_j}).
\]
That is,  $\{\widetilde{L}_\beta(w_{k_j})\}$ is bounded from below, which further implies     convergence of $\{\widetilde{L}_\beta(w_k)\}$ based on Lemma \ref{Sec3-theore1}.

Now, summing the inequality (\ref{bj-5}) over $k=0,1,\cdots,\infty$, we have by the convergence of   $\{\widetilde{L}_\beta(w_k)\}$ that
\[
\zeta_1\sum_{k=0}^{\infty}\|\triangle\mathbf{x}_{k+1}\|_G^2+
 \zeta_2 \sum_{k=0}^{\infty}\|\triangle \mathbf{y}_{k+1}\|^2\leq
 \mathcal{L}_\beta(w_0) - \widetilde{\mathcal{L}}_\beta(w_{k+1}) <\infty,
\]
which suggests $\|\triangle\mathbf{x}_{k+1}\|_G\rightarrow 0$ and  $\|\triangle\mathbf{y}_{k+1}\|\rightarrow 0$. So,      using Lemma  \ref{lem-3} and Lemma \ref{bj-1} the following holds clearly
\begin{equation} \label{bj-9-01}
\|\triangle\lambda_{k+1}\|\leq\frac{1}{\sqrt{\sigma_B}} \|B\tr \triangle\lambda_{k+1}\|\leq  \frac{L_{g}}{\sqrt{\sigma_{B}}} \|\triangle\mathbf{y}_{k+1}\|\rightarrow 0.
\end{equation}
This completes the proof. $\ \ \ \blacksquare$

Theorem \ref{lem-34} illustrates that the augmented Lagrange function of the problem (\ref{Sec1-Prob}) is convergent, and the primal and dual residuals   converge to zero. In what follows, we would present a key theorem   about    pointwise iteration-complexity of the proposed algorithm w.r.t. the primal-dual residuals. Actually, the following first assertion implies that any accumulation point of $\{w_k\}$ is a stationary point of $\{L_\beta(w_k)\}$ compared to Definition \ref{def-1}.

\begin{theorem} \label{the0-1}
Let
 $\{w_k\}$ be generated by Algorithm \ref{algo1}. Then, under Assumptions (A1)-(A3)
\begin{itemize}
 \item It holds
\begin{align}\label{p}
\lim_{k\rightarrow \infty}d(0,\partial L_{\beta}(w^{k+1}))=0.
\end{align}
%\item   It holds
%\[
%\lim_{k\rightarrow\infty}\nabla_{\lambda}L_\beta(w_{k+1})=0,~~
%\lim_{k\rightarrow\infty}\nabla_{\mathbf{y}}L_\beta(w_{k+1})=0.
%\]
%Moreover, there exists $\overline{d}_{k+1}\in \partial_{\mathbf{x}}L_\beta(w_{k+1})$ such that $\lim\limits_{k\rightarrow\infty}\overline{d}_{k+1}=0.$
\item
The sequence $\{f(\mathbf{x}_{k+1})+g(\mathbf{y}_{k+1})\}$ is convergent.
\item
Let $C_0:=\mathcal{L}_\beta(w_0)-\underline{L}.$ Then, for any integer $k\geq 1$, there exists $j\leq k$ and $\zeta_i>0\ (i=1,2,3)$ such that
\begin{equation} \label{bj-1-120}
\|\triangle\mathbf{x}_j\|_G^2\leq  \frac{C_0}{\zeta_1 (k+1)},\quad
\|\triangle\mathbf{y}_j\|^2\leq \frac{C_0}{\zeta_2 (k+1)},\quad
\|\triangle\lambda_j\|^2\leq \frac{C_0}{\zeta_3 (k+1)}.
\end{equation}
\end{itemize}
\end{theorem}
\noindent{\bf Proof }  Using  (\ref{bj-1-02}) again, we have
\[
A\mathbf{x}_{k+1}+B\mathbf{y}_{k+1}-b=-\frac{1}{\tau+\alpha }\left(\frac{1}{\beta}\triangle\lambda_{k+1}
+B\triangle\mathbf{y}_{k+1}\right)+B\triangle\mathbf{y}_{k+1},
\]
which by the third result of Theorem \ref{lem-34} suggests
\begin{equation} \label{bj-1-020}
\lim_{k\rightarrow\infty}A\mathbf{x}_{k+1}+B\mathbf{y}_{k+1}-b=0.
\end{equation}
Therefore,
\begin{align}\label{p-lambda}
\lim_{k\rightarrow\infty}\nabla_{\lambda}L_\beta(w_{k+1})=
\lim_{k\rightarrow\infty}  -(A\mathbf{x}_{k+1}+B\mathbf{y}_{k+1}-b)=0.
\end{align}
By the first-order optimality condition of $\mathbf{y}$-subproblem, it holds
\begin{eqnarray*}
0&=& \nabla g(\mathbf{y}_{k+1})-B\tr\lambda_{k+\frac{1}{2}}+\beta B\tr\left(\mathbf{x}_{k+1}^{ad}+B\mathbf{y}_{k+1}-b\right)\\
&=& \nabla g(\mathbf{y}_{k+1})-B\tr\lambda_{k+1}+\beta B\tr\left(A\mathbf{x}_{k+1}+B\mathbf{y}_{k+1}-b\right)\\
&&+B\tr\left(\lambda_{k+1}-\lambda_{k+\frac{1}{2}}\right)+\beta B\tr(\mathbf{x}_{k+1}^{ad}-A\mathbf{x}_{k+1})\\
&=&\nabla g(\mathbf{y}_{k+1})-B\tr\lambda_{k+1}+\beta B\tr\left(A\mathbf{x}_{k+1}+B\mathbf{y}_{k+1}-b\right)\\
&&-\beta B\tr(A\mathbf{x}_{k+1}+B\mathbf{y}_{k+1}-b),
\end{eqnarray*}
 which  gives
 \begin{align}\label{p-y}
 \lim_{k\rightarrow\infty}\nabla_{\mathbf{y}}L_\beta(w_{k+1})=
\lim_{k\rightarrow\infty}  \beta B\tr(A\mathbf{x}_{k+1}+B\mathbf{y}_{k+1}-b)=0.
 \end{align}
Analogously,   by the update of $\mathbf{x}$-subproblem,   there exists $d_{k+1}\in\partial f(\mathbf{x}_{k+1}) $ such that  \begin{eqnarray*}
0&=&   d_{k+1}-A\tr \lambda_{k+1}+\beta A\tr\left(A\mathbf{x}_{k+1}+B\mathbf{y}_{k}-b\right)+G(\mathbf{x}_{k+1}-\mathbf{x}_k^{md})\\
&=&  d_{k+1}-A\tr \lambda_{k+1}+\beta A\tr\left(A\mathbf{x}_{k+1}+B\mathbf{y}_{k+1}-b\right)\\
&&+\beta A\tr B(\mathbf{y}_{k}-\mathbf{y}_{k+1}) +G(\mathbf{x}_{k+1}-\mathbf{x}_{k}-\gamma_k\triangle\mathbf{x}_{k})\\
&=&  d_{k+1}-A\tr \lambda_{k+1}+\beta A\tr\left(A\mathbf{x}_{k+1}+B\mathbf{y}_{k+1}-b\right)\\
&& -\beta A\tr B\triangle\mathbf{y}_{k+1}-G(\gamma_k\triangle\mathbf{x}_{k}-\triangle\mathbf{x}_{k+1}).
\end{eqnarray*}
By defining
\[
\overline{d}_{k+1}:= d_{k+1}-A\tr \lambda_{k+1}+\beta A\tr\left(A\mathbf{x}_{k+1}+B\mathbf{y}_{k+1}-b\right),
\]
we   have $\overline{d}_{k+1}\in \partial_{\mathbf{x}}L_\beta(w_{k+1})$ and furthermore
\begin{align}\label{p-x}
\lim_{k\rightarrow\infty}\overline{d}_{k+1}
=\lim_{k\rightarrow\infty}\left[\beta A\tr B\triangle\mathbf{y}_{k+1}+G(\gamma_k\triangle\mathbf{x}_{k}-\triangle\mathbf{x}_{k+1})\right]
=0.
\end{align}
 Thus, it follows from \eqref{p-lambda}, \eqref{p-y} and \eqref{p-x} that \eqref{p} holds.

For the second assertion,    it holds by   (\ref{bj-1-020}) that
\[
f(\mathbf{x}_{k+1})+g(\mathbf{y}_{k+1})
=L_\beta(w_{k+1}) +\langle\lambda, A\mathbf{x}_{k+1}+B\mathbf{y}_{k+1}-b\rangle +\frac{\beta}{2}\|A\mathbf{x}_{k+1}+B\mathbf{y}_{k+1}-b\|^2\rightarrow L_\beta(w_{k+1}).
\]
So, the   sequence $\{f(\mathbf{x}_{k+1})+g(\mathbf{y}_{k+1})\}$ is     convergent by   the first conclusion of Theorem \ref{lem-34}.

We finally prove the pointwise iteration complexity in (\ref{bj-1-120}).
Using (\ref{lema31-01}) again, we have
\[
-\widetilde{L}_\beta(w_{k+1})\leq -
 \underline{L}   -\left(\frac{1}{2L_g}
-\frac{1}{2\beta\sigma_B}\right)\|B\tr\lambda_{k+1}\|^2 -\frac{\beta}{2}\left\|A\mathbf{x}_{k+1}+B\mathbf{y}_{k+1}-b
-\frac{\lambda_{k+1}}{\beta}\right\|^2\leq -\underline{L}.
\]
So, for any $k\geq 0,$ it follows from Lemma \ref{Sec3-theore1} that
\[\sum\limits_{j=0}^{k}\left(\zeta_1\|\triangle\mathbf{x}_{j}\|_G^2
+\zeta_2\|\triangle\mathbf{y}_{j}\|^2 \right)\leq  \mathcal{L}_\beta(w_0)+ \zeta_0\|\triangle\mathbf{x}_{0}\|_G^2-\underline{L}=C_0,\]
which shows
\[
\|\triangle\mathbf{x}_j\|_G^2\leq  \frac{C_0}{\zeta_1 (k+1)} \quad \textrm{and} \quad
\|\triangle\mathbf{y}_j\|^2\leq  \frac{C_0}{\zeta_2 (k+1)}.
\]
The final convergence rate bound in (\ref{bj-1-120}) can be also verified by   (\ref{bj-9-01}) with  $\zeta_3=L_g^2/\zeta_2. \ \ \ \blacksquare$

In order to reduce error   bounds of the primal-dual residuals,
the following remark provides an adaptive way  to update the parameter $\gamma_k$ related to $\zeta_1$ by making use of the so-called Nesterov's acceleration (proposed originally in \cite{Nesterov83}), and it also suggests how to choose   reasonable values of the parameters $\tau$ and $\alpha$.
\begin{remark}
By the above convergence analysis,
 if $G\succ 0$, then convergence of Algorithm \ref{algo1} can be guaranteed by   $\gamma_k\in[0,1/2)$. In such case we can update $\gamma_k$ adaptively by the following
\begin{equation} \label{jb--2}
\gamma_k=\frac{\theta_{k-1}-1}{2\theta_k},~~ \textrm{ where} ~ \theta_k=\frac{1+\sqrt{1+4\theta_{k-1}^2}}{2}~ \textrm{ with} ~ \theta_{-1}:=1.
\end{equation}
Note that   $\zeta_2=-\beta\sigma_B +\frac{1}{\tau+\alpha}[\beta\sigma_B-\frac{L_g^2}{\beta\sigma_B}]$  is inversely proportional to $(\tau+\alpha)$ since $L_g<\beta\sigma_B$. This together with the connection $\zeta_3=L_g^2/\zeta_2$ imply that we could choose $(\tau+\alpha)\rightarrow 1$ to get smaller error bound of $\|\triangle\lambda_j\|^2$ in (\ref{bj-1-120}). In the next section, related numerical experiments will show how to determine   reasonable values of $\tau$ and $\alpha$ in detail.
\end{remark}
%============================================================================
\section{Numerical Experiments}\label{numex}
%============================================================================
In this section, we apply the proposed algorithm to solve a class  of  practical examples from signal processing to investigate its numerical  performance. All experiments are performed  by using Windows 10 system and   MATLAB R2018a (64-bit) with an Intel Core i7-8700K CPU (3.70 GHz)  and 16GB memory.

Applying   Algorithm \ref{algo1} to solve (\ref{example-2}), we have by (\ref{bj-40}) that
\[
\mathbf{x}_{k+1}=\textrm{Prox}_{\|\mathbf{x} \|_{1/2}^{1/2},\sigma/\mu}\left(\mathbf{x}_k^{md}-\frac{  \beta A\tr(A\mathbf{x}_k^{md} -\mathbf{y}_k-b)-A\tr\lambda_k  }{\sigma}\right),
\]
which is the   half shrinkage operator \cite{xcxz12} defined as $\textrm{Prox}_{\|\mathbf{x} \|_{1/2}^{1/2}, \nu}(\mathbf{x})=(l_\nu(\mathbf{x}_1),l_\nu(\mathbf{x}_2),\cdots,l_\nu(\mathbf{x}_m))\tr$ where
\[
l_\nu(\mathbf{x}_i)=\left \{\begin{array}{llll}
\frac{2\mathbf{x}_i}{3}\left[1+\cos\frac{2}{3}(\pi-\phi(\mathbf{x}_i))\right], && \textrm{if }|\mathbf{x}_i|>\frac{3\sqrt[3]{2}}{4}\nu^{2/3}
,\\
0, && \textrm{otherwise},\\
\end{array}\right.
\]
and $\phi(\mathbf{x}_i)=\arccos(\frac{\nu}{8}(\frac{|\mathbf{x}_i|}{3})^{-3/2}).$ Besides, it is easy to obtain
$
\mathbf{y}_{k+1}=(c+\beta\mathbf{x}_{k+1}^{ad}-\lambda_{k+\frac{1}{2}})/(1+\beta).
$

With the purpose of  fast convergence  and making performance of   Algorithm \ref{algo1} less independent on an  initial guess of the penalty parameter $\beta$, as suggested by He et al.\cite{HeYangWang2000} we would adopt  the following  technique to update it adaptively:
\begin{equation}  \label{Sec32-aa}
\beta_{k+1}=\left \{\begin{array}{lllllll}
\eta^{\textrm{incr}}\beta_k& &&&& \textrm{if}\ \|r_k\|_2> \nu \|s_k\|_2,\\
{\beta_k}/{\eta^{\textrm{decr}}}& &&&& \textrm{if}\ \|s_k\|_2> \nu \|r_k\|_2,\\
\beta_k& &&&& \textrm{otherwise},\\
\end{array}\right.
\end{equation}
where $\nu, \eta^{\textrm{incr}}$ and $ \eta^{\textrm{decr}}$ are three positive parameters
with suggested values larger than $1$, for instance, $\nu=10, \eta^{\textrm{incr}}=\eta^{\textrm{decr}}=2$.
For Algorithm \ref{algo1} to solve   (\ref{Sec1-Prob}) we have
\begin{equation}  \label{Sec32-abja}
\|r_k\|= \left\|A\mathbf{x}_{k+1}+B\mathbf{y}_{k+1}-b \right\|
\end{equation}
and
\[
\|s_k\|= \left\|A\tr\triangle\lambda_{k+1} +\beta A\tr(A\mathbf{x}_{k+1}+B\mathbf{y}_{k}-b) +G(\triangle\mathbf{x}_{k+1}-\gamma_k\triangle\mathbf{x}_{k})\right\|,
\]
which   represent  the   equality constrained  error and the optimality error, respectively. Here, it is easy to check that $0\in \partial f(\mathbf{x}_{k+1})-A\tr \lambda_{k+1}+s_k.$  In order to satisfy Assumption (A2),   we  need to update   $\beta=\min\left\{\beta_{k+1}, {\frac{1.01L_g}{\sqrt{1-\tau-\alpha}\sigma_B}}\right\}$ at each iteration. As for the problem (\ref{example-2}), we have $L_g=1$ and $ \sigma_B=1$.  If not specified, the initial penalty parameter $\beta_0$ is chosen as $0.04$,   the starting points $(\mathbf{x}_0,\mathbf{y}_0)$ and $\lambda_0$ are respectively set as zero and ones vector  with proper dimensions, and  the matrix  $G=\sigma I -\beta A\tr A$ with $\sigma=1.01 \beta\|A\tr A\|$.  The parameter $\gamma_k$ is updated adaptively according to (\ref{jb--2}). Throughout  we use the following  stopping criterion as mentioned in   \cite{LiPONG15}    to terminate Algorithm \ref{algo1}:
\begin{equation} \label{stop}
\textrm{IRE(k)}:= \frac{\max\left\{\|\mathbf{x}_k-\mathbf{x}_{k-1}\|, \|\mathbf{y}_k-\mathbf{y}_{k-1}\|,\|\lambda_k-\lambda_{k-1}\|\right\}}{\max\left\{\|\mathbf{x}_{k-1}\|, \|\mathbf{y}_{k-1}\|,\|\lambda_{k-1}\|,1\right\}}<\epsilon,
\end{equation}
where $\epsilon$ is a given tolerance error.  Note that this stopping criterion corresponds to the pointwise iteration complexity shown in (\ref{bj-1-120}), so such stopping criterion is well defined.

As the first experiment, we consider the reformulated sparse signal recovery problem (\ref{example-2}) with an original signal $x\in \mathcal{R}^{3072}$   containing 160 spikes with amplitude $\pm 1$. The measurement matrix $A\in \mathcal{R}^{1024\times 3072}$ is drawn  firstly from  the standard
norm distribution $\mathcal{N}(0,1)$ and then each of its column is normalized. Specifically, we use the following MATLAB codes to generate the original signal $x_{\textrm{orig}}$, the data $A,  c$ and $\mu$:
\[\begin{array}{l}
\verb"randn('state', 0); rand('state',0);"\\
\verb"l = 1024; m = 3072;"\\
\verb"T = 160;  % number of spikes"\\
\verb"x_orig = zeros(m,1); q = randperm(m);"\\
\verb"x_orig(q(1:T)) = sign(randn(T,1));  % original signal" \\
\verb"A = randn(l,m);"\\
\verb"A = A*spdiags(1./sqrt(sum(A.^2))',0,m,m);  % normalize columns"\\
\verb"sig = 0.01;  % noise standard deviation"\\
\verb"c = A*x_orig + sig*randn(l,1);  % noisy observations"\\
\verb"mu_max = norm( A'*c,'inf');"\\
\verb"mu = 0.1*mu_max;  % regularization parameter"\\
\end{array}
\]
Under   tolerance $\epsilon=10^{-15}$, we   test the effect of   parameters $(\tau,\alpha)$   restricted in (\ref{al-111}) on  the  numerical performance of Algorithm \ref{algo1} (In fact, we choose parameter values around   $(\tau,\alpha)=(0.3,0.32)$, because we find it performs  slightly better than some pairs after running a lot of  values restricted in (\ref{al-111}) by the aid of two level for loops in MATLAB). We also randomly choose four pairs of $(\tau,\alpha)$ to carry out related experiments.

Table 1 reports some  computational results of several quality measurements, including ``IT'', ``CPU'', ``IRE'', ``EQU''  which  denote respectively the   iteration number, the CPU time in seconds, the final relative iterative error IRE(k) defined in (\ref{stop}) and the final feasibility error $\|r_k\|$ defined in (\ref{Sec32-abja}).  We use $l_2$-error (defined as $\frac{\|x_k-x_{\textrm{orig}}\|}{\|x_{\textrm{orig}}\|}$) to represent   the  relative error   to measure   recovery quality of a signal. As shown in Table 1, setting $(\tau,\alpha)=(0.65, 0.32)$ would be  a reasonable choice for Algorithm \ref{algo1} to solve the  problem (\ref{example-2}),  because  in such  a choice the iteration number and the CPU time are relatively smaller  while  reported  results in  each of the last three columns  are nearly the same when   the stopping criterion is satisfied. Hence, in the following experiments, we   use Algorithm \ref{algo1} with default parameters $(\tau,\alpha)=(0.65, 0.32)$.

{\small
\begin{center}
\begin{tabular}{cccccc}
\hline
 $(\tau,\alpha)$  & IT  & CPU  & IRE &    EQU  & $l_2$-error\\
\hline
(0.3, 0.10)& 472 & 22.67 & 9.47e-16&6.10e-14 & 6.79e-2\\
(0.3, 0.15)& 443 & 21.12 & 9.23e-16&5.18e-14 & 6.79e-2 \\
(0.3, 0.20)& 500 & 24.11 & 9.77e-16&4.82e-14  &6.79e-2 \\
(0.3,  0.25)&379 & 18.23 & 9.79e-16&6.33e-14 &  6.79e-2\\
(0.3,  0.30)& 350  & 16.79 & 9.78e-16&5.88e-14 & 6.79e-2 \\
(0.3,  0.32)& 344 & 16.66  & 9.45e-16 &6.33e-14 & 6.79e-2 \\
(0.3,  0.35)& 544 & 26.17  & 9.79e-16 &3.28e-14  & 6.79e-2\\
(0.3,  0.40)& 510 & 24.30  & 8.87e-16 &5.11e-14  & 6.79e-2\\
(0.3,  0.45)& 485 & 23.28  & 9.93e-16 &3.34e-14  & 6.79e-2\\
(0.3,  0.50)& 458 & 21.95  & 9.34e-16 &3.24e-14  & 6.79e-2\\
(0.3,  0.55)& 433& 20.86   & 9.58e-16 &3.29e-14  & 6.79e-2\\
(0.3,  0.60)& 413 & 20.30  & 9.45e-16 &3.24e-14  & 6.79e-2\\
(0.3,  0.65)& 396 & 19.24  & 9.33e-16 &3.26e-14  & 6.79e-2\\
(0.3,  0.68)& 387 &18.66  & 8.91e-16 &3.37e-14   & 6.79e-2\\
\hline
(-0.3, 0.32)&$-$  & $-$  & $-$& $-$  &$-$ \\
(-0.2, 0.32)&$-$  & $-$  & $-$& $-$ & $-$\\
(-0.1, 0.32)&695  &33.51& 9.95e-16& 5.41e-14 & 6.79e-2\\
(0, 0.32)  &498 &23.77  & 9.99e-16& 8.38e-14 & 6.79e-2\\
(0.1, 0.32)&697	  & 33.49&9.33e-16& 5.68e-14 & 6.79e-2\\
(0.2, 0.32)&391   & 18.93  & 9.54e-16& 7.85e-14 & 6.79e-2 \\
(0.3, 0.32)&  344  & 16.62 & 9.45e-16 &6.33e-14 & 6.79e-2\\
(0.4, 0.32)&502   & 24.16  & 8.85e-16& 4.01e-14 & 6.79e-2 \\
(0.5, 0.32)&451   &21.61   & 9.49e-16& 2.97e-14  & 6.79e-2\\
(0.6, 0.32)&404   & 19.52  & 9.64e-16& 3.74e-14 & 6.79e-2\\
(0.62, 0.32)&397	 &19.20   & 9.71e-16 & 4.34e-14 & 6.79e-2 \\
(0.65, 0.32)& \textbf{279}& \textbf{13.52}& 9.98e-16 &3.19e-14  & 6.79e-2\\
(0.67, 0.32)&318	 &15.59  &  8.43e-16& 3.11e-14  &6.79e-2 \\
\hline
(0.90, 0.05)&396	 &19.06 & 9.26e-16& 3.27e-14  &6.79e-2 \\
(0.80, 0.15)&396 	 & 19.02&  9.33e-16&  3.28e-14  & 6.79e-2\\
(0.01, 0.90)&411 	 & 19.78&  9.63e-16&  3.10e-14  &6.79e-2 \\
(0.05, 0.70)& 485	& 24.62&  9.62e-16&  3.10e-14  &  6.79e-2\\
\hline
\end{tabular}
\end{center}}
\begin{center}
Table 1:\ Results\footnote{``$-$'' means that the stopping criterion is not satisfied after 800 iterations, and the bold number   in that row indicate the best results obtained  by changing $(\tau,\alpha)$ belong to $(0,1)$.} of  Algorithm \ref{algo1}  with  different $(\tau,\alpha)$ for solving  problem (\ref{example-2}).
\end{center}

\begin{figure}[htbp]
 \begin{minipage}{1\textwidth}
 \def\figurename{\footnotesize Fig.}
 \centering
\resizebox{16cm}{5.2cm}{\includegraphics{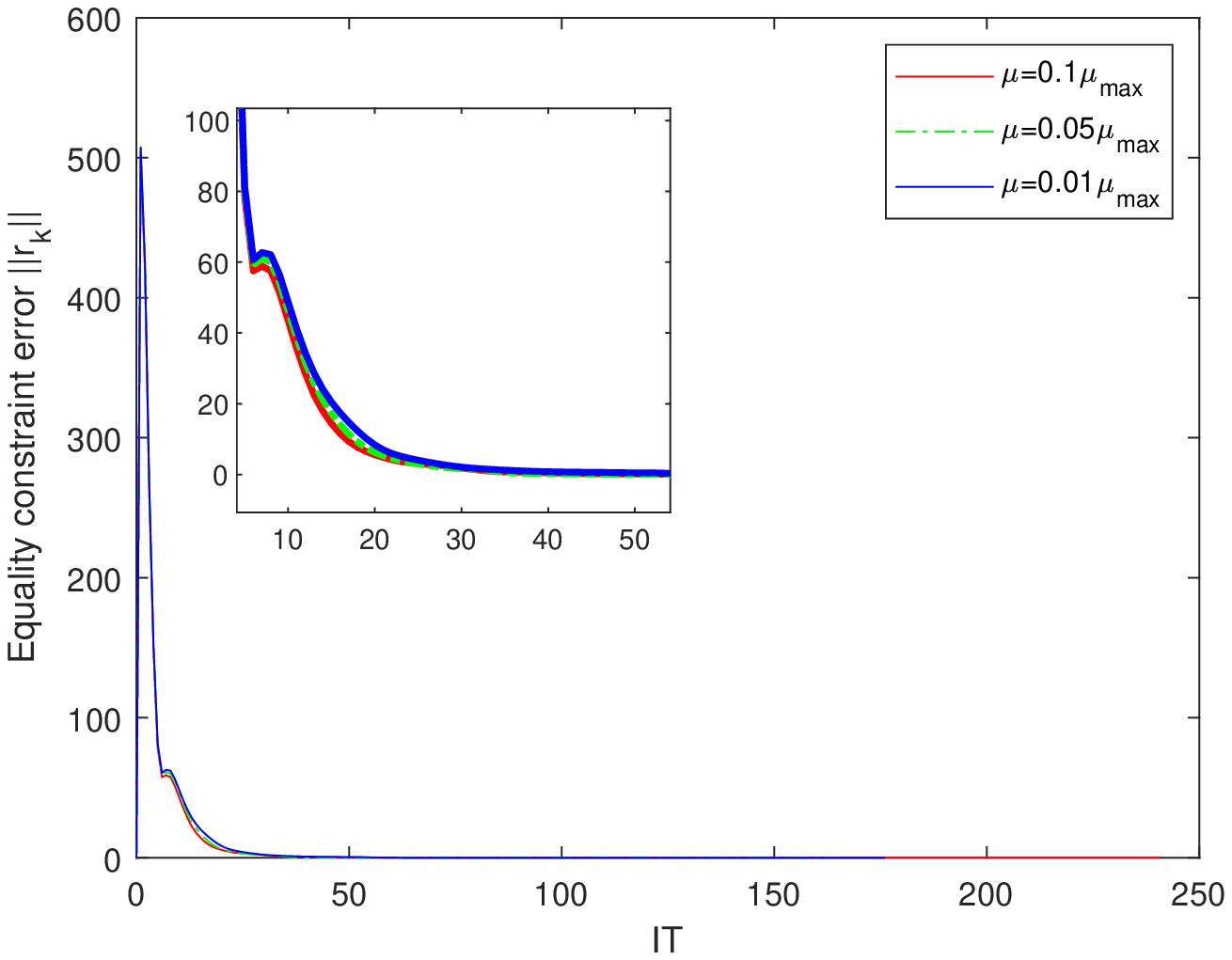}
\includegraphics{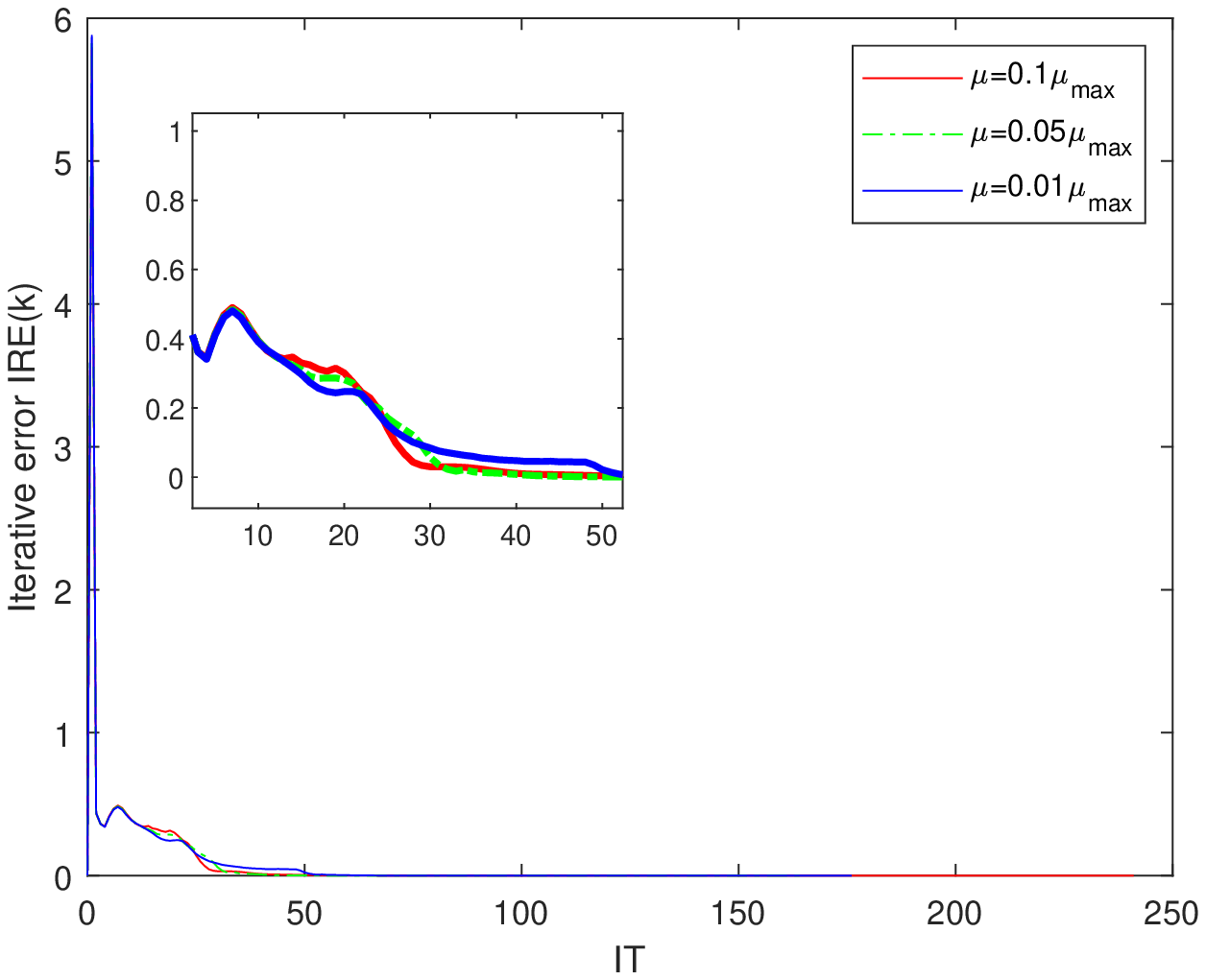}\includegraphics{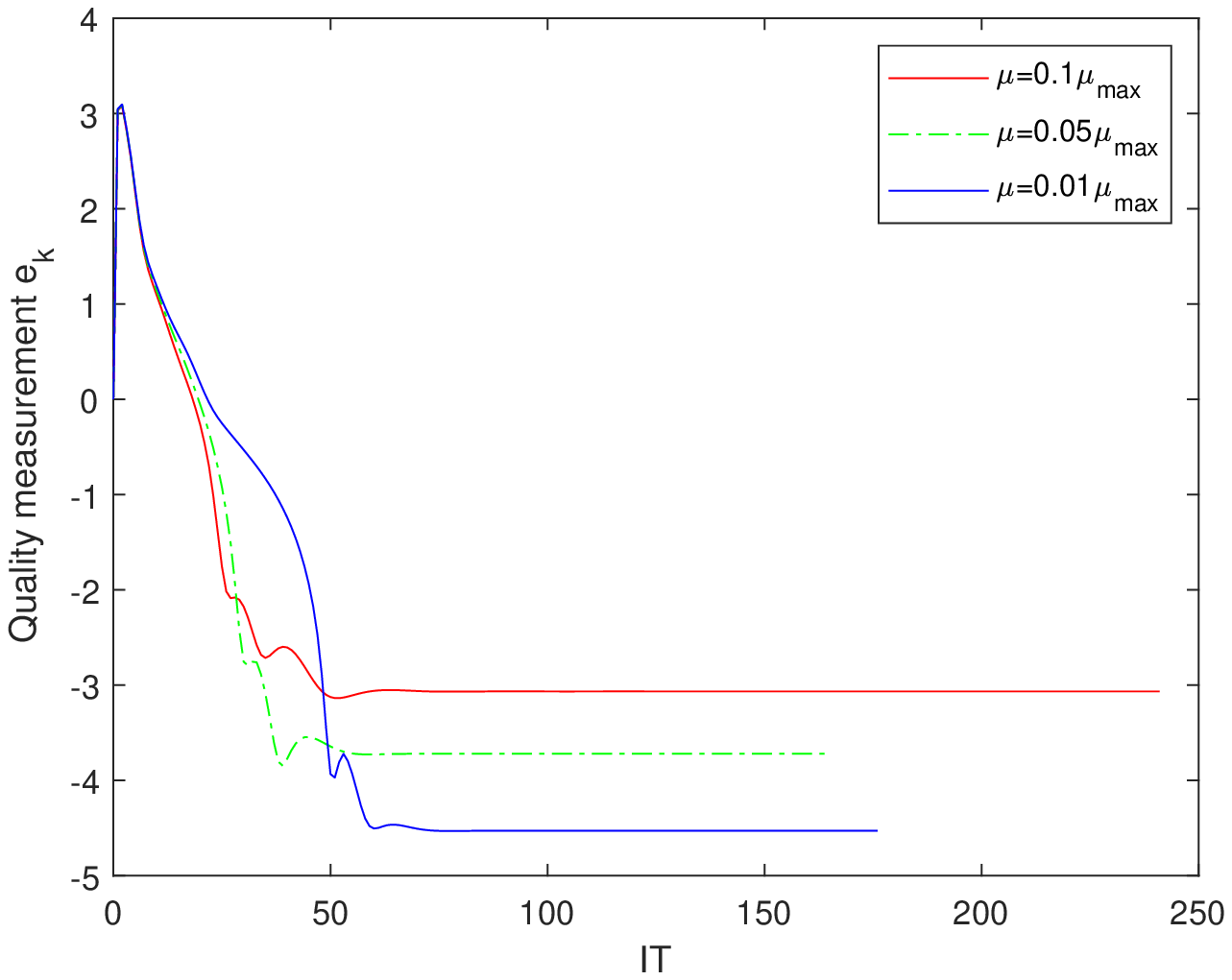}}
\caption{\footnotesize Convergence tendency of the equality constraint error $\|r_k\|$ (left), the iterative error IRE(k) (middle) and the recovery signal   quality $e_k$ (right)  by   Algorithm \ref{algo1} for solving  the problem (\ref{example-2}) with   $(l,m)=(2048,5000)$ but with  different regularization factors. }
   \end{minipage}
\end{figure}

Next, we use the aforementioned codes to investigate the effect of regularization parameter $\mu$ on Algorithm \ref{algo1}  for solving the problem (\ref{example-2}) with a large   data  $A\in \mathcal{R}^{2048\times 5000}$ and the same spikes, but the tolerance is set as $\epsilon=10^{-12}$. Fig. 1 depicts   convergence behaviors of  the equality constraint error $\|r_k\|$, the iterative error IRE(k)  and the recovery signal quality  $e_k:=\log_{10}\frac{\|\mathbf{x}_k-\mathbf{x}_{\textrm{orig}}\|}{\|\mathbf{x}_{\textrm{orig}}\|}$ along the iteration process after applying Algorithm \ref{algo1} with $\mu=0.1\mu_{\max},\  0.05\mu_{\max},\ 0.01\mu_{\max}$, respectively. Fig. 2 also presents the results to visualize the recovery quality of  the  signal versus  the original signal, where the upper-left plot shows the minimum energy reconstruction signal $A^{\dagger}c$ (which is the point satisfying $A\tr A\mathbf{x}=A\tr c$) versus the original signal. An outstanding observation from Fig. 1 is that the smaller the value of $\mu$ is, the smaller the iteration number is (and the better   the recovery quality of the signal is). After identifying the nonzero positions in the reconstructed signal,  it  always has the correct number
of spikes for the case with $\mu=0.01\mu_{\max}$  and is   closer to the original noiseless signal.

\begin{figure}[htbp]
 \begin{minipage}{1\textwidth}
 \def\figurename{\footnotesize Fig.}
 \centering
\resizebox{16.5cm}{9.1cm}{\includegraphics{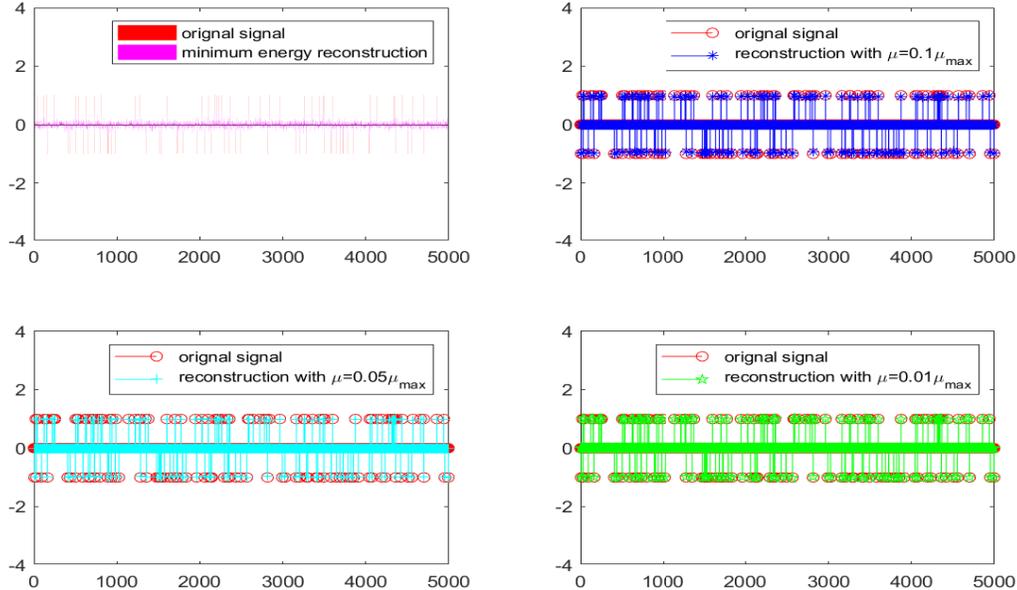}}
\caption{\footnotesize Comparison between the original signal and reconstructed signal   by   Algorithm \ref{algo1} for solving   problem (\ref{example-2}) with   $(l,m)=(2048,5000)$ but with different regularization factor $\mu$. }
   \end{minipage}
\end{figure}

{\small
\begin{center}
\begin{tabular}{c|cccc|cccc}
\hline
   &    &   $l_{1/2}$   &regularizer&   & & $l_1$  & regularizer&\\
\hline
 $(l,m)$  & IT  & CPU &   EQU  & $l_2$ error  &IT & CPU & EQU  & $l_2$ error\\
\hline
(1024, 3000)& {\bf 358} &{\bf 16.37} &4.60e-14 &{\bf 1.20e-2} & 501& 22.60&1.93e-14&3.70e-2 \\
(1024, 4000)& {\bf 367}  & {\bf 25.93}  &  4.18e-14 & {\bf 1.28e-2}  &507 & 35.611&4.78e-14& 4.26e-2 \\
&   &   &    &   & & & & \\

(2048, 5000)& {\bf  215} & {\bf 30.65}  &  2.84e-14 & {\bf 1.08e-2}   & 250&36.01 &3.27e-14&2.66e-2 \\
(2048, 6000)& {\bf  222} & {\bf 41.37}  &  3.98e-14 &{\bf  1.20e-2}  & 266&49.59&3.47e-14& 3.07e-2 \\
&   &   &    &   & & & & \\

(3000, 7000)&{\bf  201}	 &{\bf 58.08}    &  2.71e-14 & {\bf 1.17e-2}   &231  &  66.91 & 2.93e-14& 2.60e-2\\
(3000, 8000)& {\bf 205}  &{\bf 71.36}&  3.77e-14 & {\bf 1.10e-2}  & 230&79.55 & 3.62e-14& 2.58e-2 \\
&   &   &    &   & & & & \\

(4000, 9000)& {\bf 199}	 &{\bf  97.52}  &  3.29e-14 & {\bf  1.11e-2}   &  231 &112.94  &  2.87e-14&  2.69e-2\\
(4000, 10000)&  {\bf  202}   &{\bf 118.62} & 2.37e-14 &  {\bf 1.03e-2}  & 231&135.72  &  2.91e-14&  2.51e-2 \\
\hline
\end{tabular}
\end{center}}
\begin{center}
Table 2:\ Results  of  Algorithm \ref{algo1}  for   (\ref{example-2}) with different regularization terms and   dimensions.
\end{center}

In the following, we use the proposed algorithm to solve two different cases of the sparse signal recovery problem to investigate which regularization term performs better:
 \textbf{Case (i)} the  convex problem (\ref{example-1})\footnote{Note that this is also a special case of (\ref{Sec1-Prob}) with $f(\mathbf{x})=\mu\|\mathbf{x} \|_1,
  g(\mathbf{y})=\frac{1}{2}\|\mathbf{y}-c\|^2, B=-I$ and $ b=0.$} with $l_1$ regularization term;
  \textbf{Case (ii)} the nonconvex problem   (\ref{example-2}) with $l_{1/2}$ regularization term.  Table 2 reports some numerical results, where the problem dimension comes from 3000 to 10000 w.r.t the dimension of the signal,   the regularization parameter is fixed as $\mu=0.01\mu_{\max}$ and Algorithm \ref{algo1} is terminated  under   tolerance $\epsilon=10^{-15}$ with maximal iteration numbers 1000. Fig. 3 depicts comparison results between the original signal and the reconstructed signal for the signal dimension $m=10000$. First of all, it can be seen from results   in Table 2 that the proposed algorithm is feasible for solving both the nonconvex and convex sparse signal recovery problem, especially for the large-scale problem. Besides, an   obvious observation from Table 2 is that using $l_{1/2}$ regularizer is significantly better than $l_1$ regularizer to recover a signal, which could be checked from reported results of the iteration number, the CPU time and the recovery quality (i.e., $l_2$ error).

 \begin{figure}[htbp]
 \begin{minipage}{1\textwidth}
 \def\figurename{\footnotesize Fig.}
 \centering
\resizebox{14.5cm}{14cm}{\includegraphics{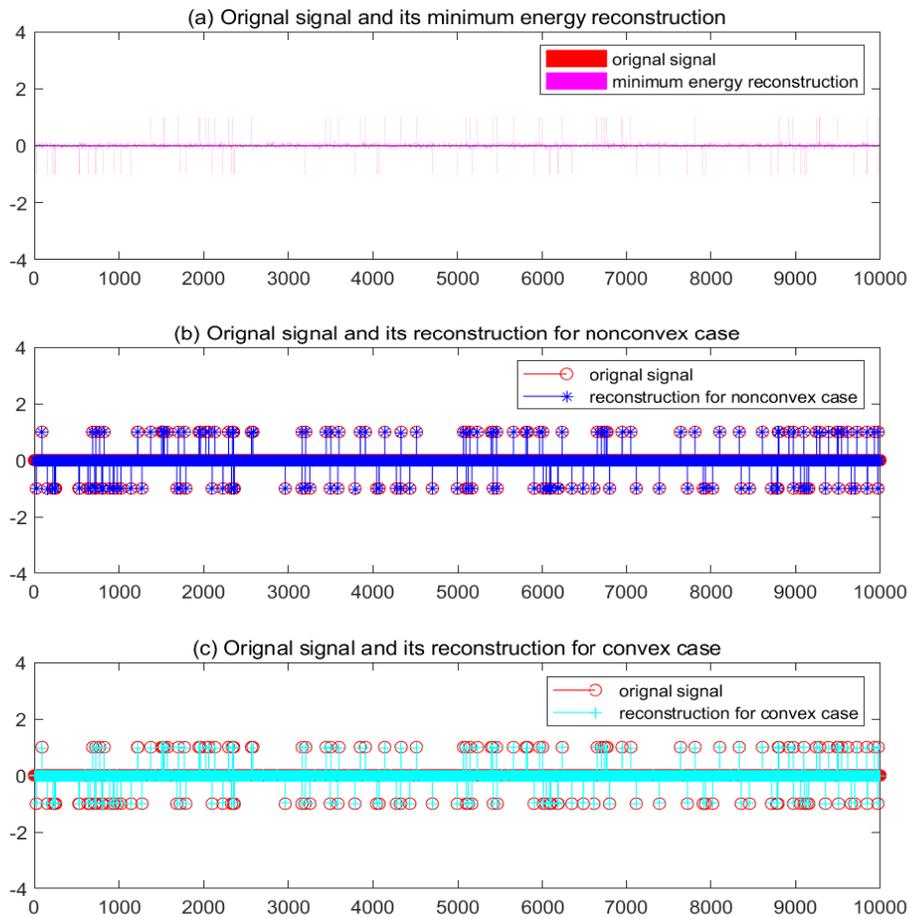}}
\caption{\footnotesize Original signal and reconstructed signal   by   Algorithm \ref{algo1}  for solving the sparse signal recovery problem with   $(l,m)=(4000,10000)$ but with different regularization terms. }
   \end{minipage}
\end{figure}

Finally, we would apply the proposed algorithm to solve the direction-of-arrival (DOA) estimation problem   \cite{MCIE05} with a single snapshot.
Here we consider a uniformly  linear array of $M=100$ sensors  with half-wavelength  elements spacing.
Let $\boldsymbol{\theta}=[\theta_1,\cdots,\theta_L]\tr$ denote the $L$ angles of interest in $[-\frac{\pi}{2}, \frac{\pi}{2}]$.
Denote $\mathbf{x}=[x_1,\cdots,x_L]\tr$ as the amplitudes of the potential signals from the $L$ incoming angles. Thus, the received signal   at the sensor array is given by:
$\mathbf{y}=A\mathbf{x}+\boldsymbol{n}$,
where
$\mathbf{y}=[y_1,\cdots,y_M]\tr$,
$\boldsymbol{n}=[n_1,\cdots,n_M]\tr$,
the steering matrix
$A=[\mathbf{a}(\theta_1),\cdots,\mathbf{a}(\theta_L)]$ and
$\mathbf{a}(\theta_l)=[1,
\exp(-j\pi \sin(\theta_l)),
\cdots,
\exp(-j\pi (M-1)\sin(\theta_l))]\tr$.%,

\begin{figure}[htbp]
 \begin{minipage}{1\textwidth}
 \def\figurename{\footnotesize Fig.}
 \centering
\resizebox{13cm}{5cm}{\includegraphics{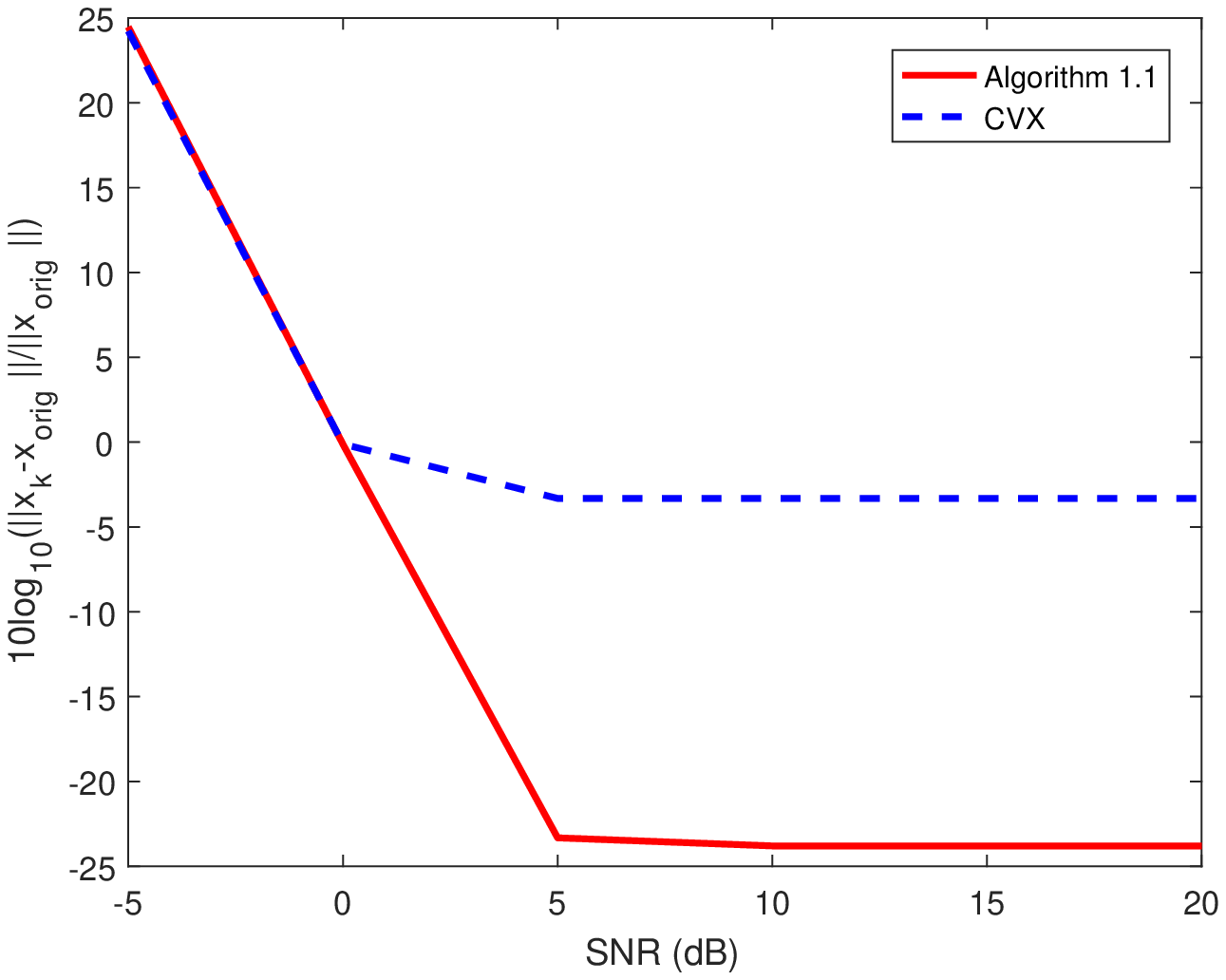}
\includegraphics{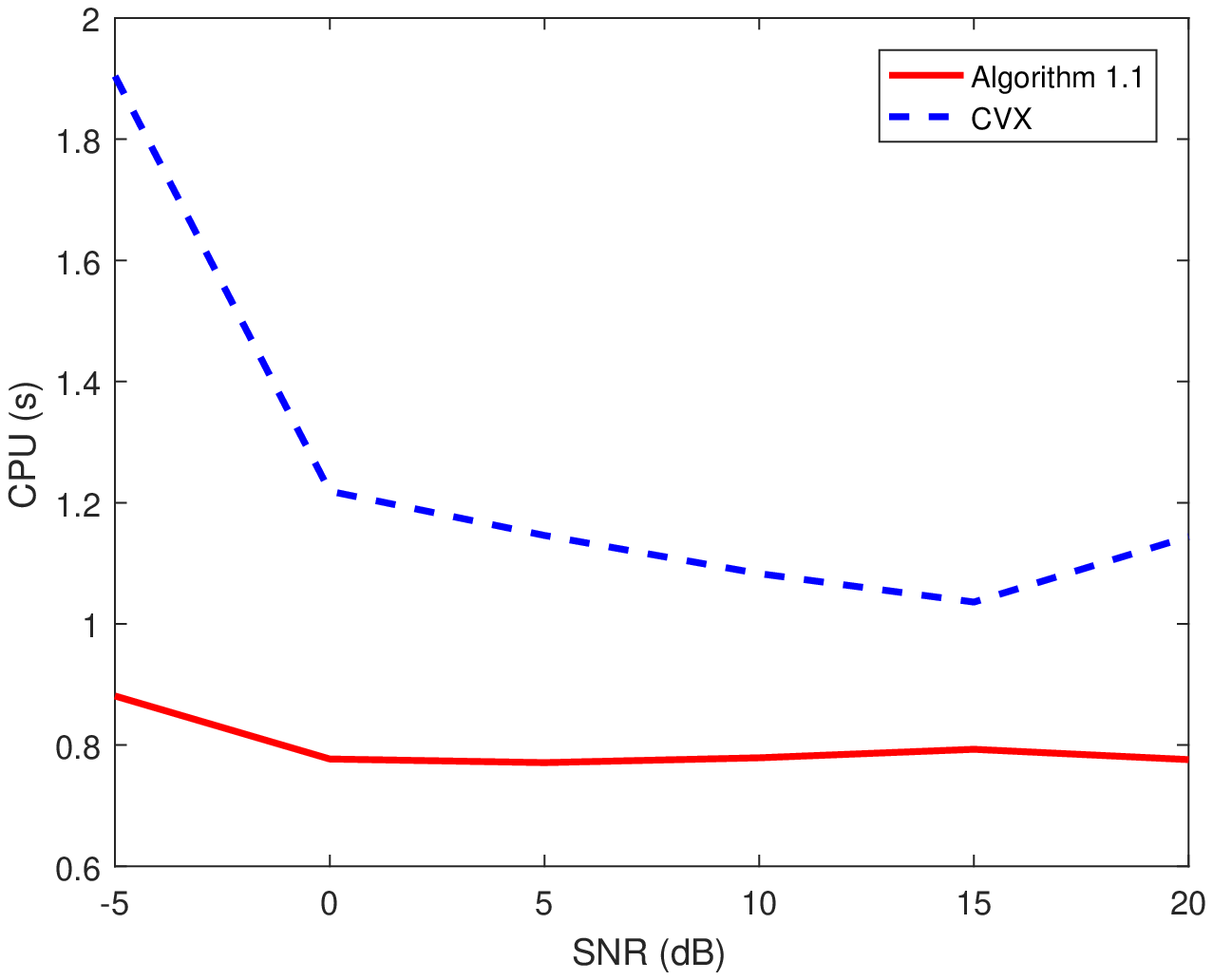}}
\caption{\footnotesize The left and right are the   errors   and   average run time versus different SNRs, respectively.}
   \end{minipage}
\end{figure}

\begin{figure}[htbp]
 \begin{minipage}{1\textwidth}
 \def\figurename{\footnotesize Fig.}
 \centering
\resizebox{14cm}{6cm}{\includegraphics{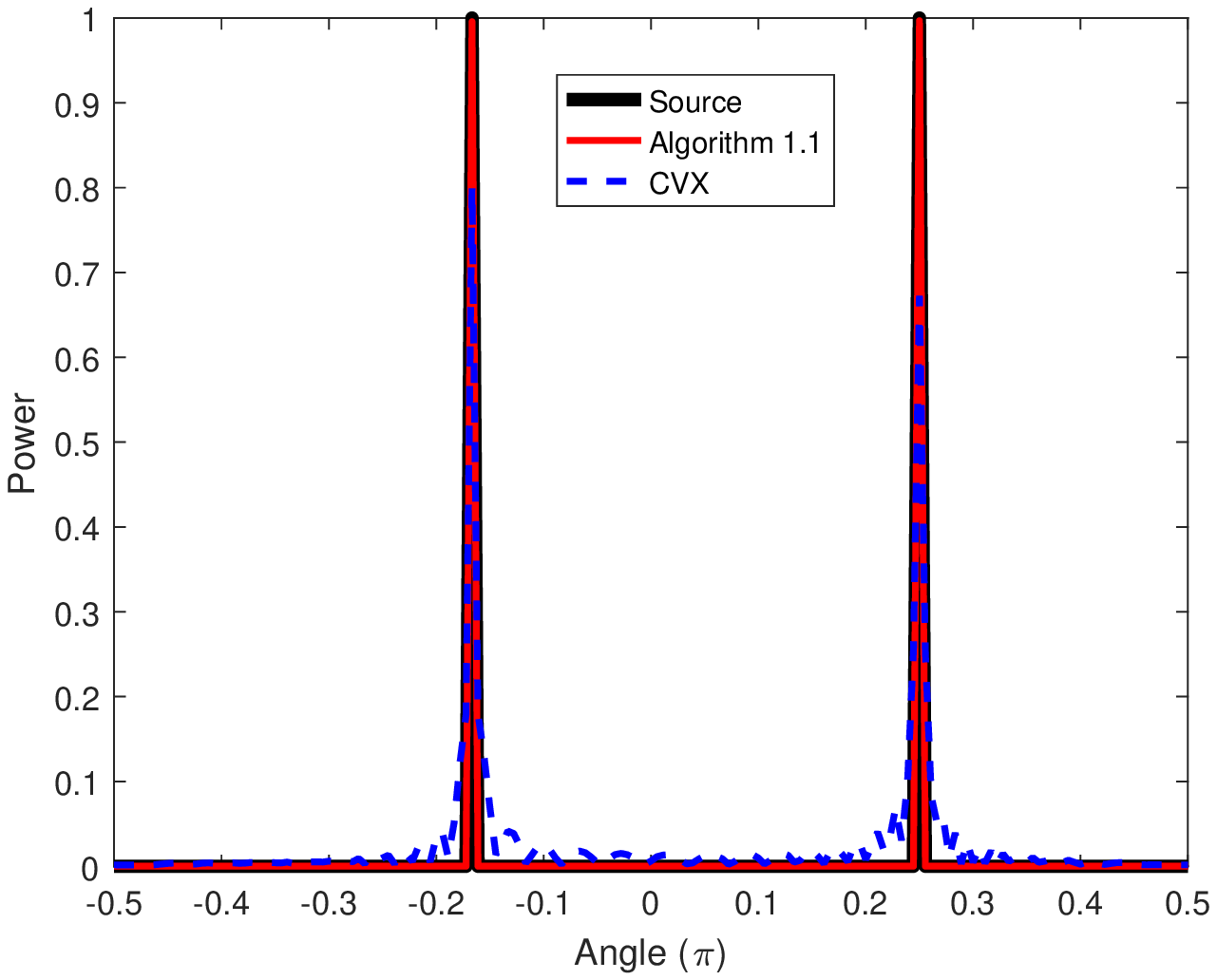}}
\caption{\footnotesize DOAs at SNR$=5$dB.}
   \end{minipage}
\end{figure}

We consider the narrowband scenario with $K=2$ uncorrelated far-field source signals   with normalized DOA parameters $-\frac{\pi}{6}$ and $\frac{\pi}{4}$. To run the proposed method, we divide the potential angle region $[-\frac{\pi}{2}, \frac{\pi}{2}]$ into $L=180$ uniformly discrete grid points, i.e., $\boldsymbol{\theta}=\frac{\pi}{180}[-90,-89,\cdots,89,90]\tr$.
When the signal-to-noise-ratio (SNR) varies from $-5$dB to $20$dB, i.e.,  $\{-5,0,5,10,15,20\}$dB, we implement the proposed method and the well-known CVX\footnote{Avaliable at: http://cvxr.com/cvx/.} toolbox for $100$ Monte Carlo runs, and compute their root mean square errors and running time, as plotted in Fig. 4. For visible comparison, we plot the result from one Monte Carlo   in the case of $5$dB, as shown in Fig. 5. From Figs. 4-5, we can see that:
\begin{itemize}
\item
The accuracy of the two methods increases with the increase of SNR;
\item
The implementation of the proposed method is faster than that of the CVX method.
\item
In terms of DOA resolution and the estimation accuracy of the incoming signal power, the proposed method is better than that of CVX.
\end{itemize}

%=============================
\section{Conclusion remarks} \label{Fin-5}
%=============================
In this paper, we construct a   symmetric alternating direction method of multipliers for solving a family of possibly nonconvex nonxmooth optimization problems. Two different acceleration techniques are  designed for fast convergence. Under proper assumptions, convergence of the proposed algorithm as well as its pointwise iteration complexity are analyzed in detail. By testing the so-called sparse signal recovery problem in signal processing  with nonconvex/convex regularization terms and by using adaptively updating strategy for the penalty parameter, a number of numerical results demonstrate the feasibility and  efficiency of the new algorithm and further show that the $l_{1/2}$ regularization term is better than the $l_{1}$ regularization term in terms of CPU time, iteration number and recovery error. Our future work will focus on  solving stochastic nonconvex optimization problems by using a similar first-order algorithm  to ADMM.

\end{document}